\documentclass[12pt]{amsart}
\usepackage{graphicx} % Required for inserting images
\usepackage{amsmath,amsthm,amsfonts,bbm,bm,color,cite,enumerate}
\usepackage[foot]{amsaddr}
\newtheorem{theorem}{Theorem}[section]
\newtheorem{cor}[theorem]{Corollary}
\newtheorem{lem}[theorem]{Lemma}

\theoremstyle{definition}

\theoremstyle{remark}

\numberwithin{equation}{section}
\usepackage[margin=1in]{geometry}
\usepackage[nameinlink,noabbrev]{cleveref}

\newtheorem{assumption}{Assumption}
\Crefname{assumption}{Assumption}{Assumptions}

\title[Feynman--Kac--Lyapunov Method]{Constructing Lyapunov functions of stochastic differential equations by solving the Feynman--Kac equation}
\author{Weiqi Chu}
\author{Ning Jiang}
\author{Yao Li}
\thanks{Department of Mathematics and Statistics, University of Massachusetts Amherst, Amherst, Massachusetts, the United States, 01003}
\email{weiqichu@umass.edu, ningjiang@umass.edu, liyao@umass.edu}

\date{\today}

\begin{document}

\begin{abstract}
The Lyapunov function, which is often constructed on an ad hoc basis, plays a central role in proving the stochastic stability of stochastic differential equations (SDEs). In this paper, we present a systematic approach for constructing the Lyapunov function for a class of ``driving–driven'' SDEs. By assuming that the Lyapunov function has a product form, we show that the prefactor function associated with the driving variables must satisfy a Feynman--Kac equation. We establish sufficient conditions under which this equation has a strictly positive solution and construct the Lyapunov function based on it. We discuss several special cases in which these conditions can be relaxed. In addition, we solve the Feynman--Kac equation numerically, integrating its probabilistic representation with a deep learning approach. We develop a neural network solver for high-dimensional Feynman--Kac equations and demonstrate its effectiveness through numerical experiments.
\end{abstract}
\maketitle

\section{Introduction}
For a stochastic differential equation (SDE), a natural question is whether it is stochastically stable, or more precisely, whether it has a unique invariant probability measure that is ergodic. Throughout this paper, we use the term \emph{stochastic stability} to refer to the existence and uniqueness of an invariant probability measure and to either geometric or power-law convergence toward this invariant measure. The stochastic stability of SDEs has been a central topic of study for decades. Despite extensive theoretical progress, proving stochastic stability remains challenging for many high-dimensional systems encountered in applications. The Lyapunov function method is the most widely used tool for such analyses, yet there is no universal procedure for constructing Lyapunov functions for arbitrary SDEs.

Consider a $d$-dimensional SDE
\begin{equation}
    \label{SDE1}
    \mathrm{d}X_t = f(X_t)\,\mathrm{d}t + \sigma(X_t)\,\mathrm{d}W_t, 
\end{equation}
where $f$ is a vector field from $\mathbb{R}^d$ to itself, $\sigma \in \mathbb{R}^{d\times n}$ is a matrix-valued function, and $W_t$ is the Wiener process on $\mathbb{R}^n$. Let $\mathcal{L}$ be the infinitesimal generator of equation \eqref{SDE1}. If a function $V: \mathbb{R}^d \rightarrow [0, \infty)$ satisfies
\begin{equation}
\mathcal{L}V \leq - a V^\alpha + b
\end{equation}
for some positive constants $a, b$, and $\alpha$, we say that $V$ is a Lyapunov function for $X_t$.

Stochastic stability typically follows from the existence of an appropriate Lyapunov function together with a minorization condition near the ``bottom'' of the Lyapunov function $\{ V \leq c\}$ for some constant $c > 0$. The minorization condition ensures that the stochastic process \eqref{SDE1} possesses a common future when initialized from a certain set.
More precisely, let $P^t(x, \cdot)$ be the transition kernel of a Markov process $X_t$. We say $X_t$ satisfies the {\it minorization condition} if there exist a set $C \subset \mathbb{R}^d$, a nontrivial probability measure $\nu$, and constants $t_0, \delta > 0$, such that 
\begin{equation}
    \label{smallset}
    P^{t_0}(x, \cdot) \geq \delta \nu(\cdot), \quad \forall x \in C.
\end{equation}
If the minorization condition \eqref{smallset} holds, then the equation \eqref{SDE1} achieves exponential convergence to a unique invariant probability measure when $\alpha = 1$, and power-law convergence when $0 < \alpha < 1$. We refer to \cite{meyn1993MC, khasminskii2012stochastic} for results on the stochastic stability of Markov processes, to \cite{Oksendal, khasminskii2012stochastic} for SDEs, and to \cite{hairer2006ergodicity, herzog2015noise1, herzog2015noise2, mattingly2002ergodicity} for recent developments on stochastic stability, including various relaxations of the minorization condition.

For an SDE, the minorization condition usually follows from accessibility, which means any open set can be visited by the stochastic process with strictly positive probability, plus the hypoellipticity or H\"ormander's condition. See, for example, Theorem 7.4.20 of \cite{stroock2008partial} for H\"ormander's condition and \cite{hairer2006ergodicity} for the minorization condition.
For a given system, it is usually not difficult to check its minorization condition. In contrast, there is no systematic approach for constructing a Lyapunov function. In general, its construction depends on the understanding of the dynamics of the SDE; intuitively, the Lyapunov function should assign higher values to states that the process $X_t$ is less likely to visit. Some recent studies address this challenge by solving the Poisson equation to bridge the gaps (i.e., regions without properly defined Lyapunov functions) between regions with well-defined Lyapunov functions \cite{athreya2012propagating}, but this approach is not applicable to all problems.

In this paper, we provide a systematic framework, called the Feynman--Kac--Lyapunov method, for constructing Lyapunov functions for a class of high-dimensional SDEs that comprise a ``driving'' component $Y_t$ and a ``driven'' component $X_t$, where the stability of the driven component depends on the state of the driving component. Many stochastic systems in ecology, epidemiology, and statistical physics belong to this class \cite{hening2018coexistence, bedrossian2024stationary}. The one-dimensional version of this method was originally developed in \cite{hani2025nonequilibrium}, which proves the stochastic stability of a highly challenging energy cascade model. In this paper, we adopt a similar approach and extend it to more general settings. 
In the literature, the Feynman--Kac--Lyapunov framework typically involves estimating Lyapunov exponents for two-point processes \cite{baxendale1988large, bedrossian2022almost, bedrossian2025non, dolgopyat2004sample}, which requires solving for the eigenfunction of a Feynman--Kac operator. In contrast, our framework mainly emphasizes the solution of the Feynman--Kac equation on a bounded domain, employing an exterior Lyapunov function to control the unbounded region.

We illustrate the main idea of this paper using the following toy example, in which the driven component $X_t$ and the driving component $Y_t$ are scalar, and $X_t$ has a linear driving force without noise. That is
\begin{equation}  \label{toyexample}
\begin{aligned}
        \mathrm{d}X_t & = X_tc(Y_t)\,\mathrm{d}t, \\
        \mathrm{d}Y_t & = g(Y_t)\,\mathrm{d}t + \sigma\,\mathrm{d}B_t. 
\end{aligned}
\end{equation}
For the sake of simplicity, we focus on constructing a natural Lyapunov function $V$ for $X_t$. If $c(y) < 0$ for all $y$, then $|x|^\alpha$ is a Lyapunov function for all $\alpha > 0$. 
However, in general, the sign of $c$ cannot be guaranteed, and the stochastic stability of $X_t$ depends on the dynamics of $Y_t$. If, on average, $c(Y_t)$ is negative, then over a long time the expected value of $|X_t|^\alpha$ should still decrease. An intuitive approach is to introduce a prefactor $U(y)$ for the natural Lyapunov function $|x|^\alpha$, where the prefactor penalizes the region in which $c(y)$ is positive. 
This motivates a test function
\begin{equation}\label{eq: W}
W(x,y)=|x|^{\alpha}U(y),\end{equation}
where $U(y)$ is the prefactor we construct based on the dynamics of $Y_t$. 
Let $\mathcal{L}$ be the infinitesimal generator of \eqref{toyexample}. From It\^{o}'s formula, we obtain that 
\begin{equation}\label{eq: LWXY}
\begin{aligned}
    \mathcal{L}W(x,y) = |x|^{\alpha}\left[U'(y)g(y)+\frac{1}{2}\sigma^2U''(y)+\alpha c(y)U(y)\right],
\end{aligned}
\end{equation}
using the fact that $|x|=\mbox{sign}(x)x$. Since the dynamics of $Y_t$ do not depend on $X_t$, we have the infinitesimal generator of $Y_t$, which is
\begin{equation}\mathcal{L}^Y U = (U' g)+\frac{1}{2}\sigma^2U''.\end{equation}
From \eqref{eq: W} and \eqref{eq: LWXY}, we know that $W$ is a Lyapunov function of $(X_t, Y_t)$ if $U$ is a strictly positive solution to an equation as follows 
\begin{equation} \label{eq: lyapunov_1d}
\mathcal{L}^YU+\alpha cU=-q
\end{equation}
with $q/U > a > 0$ uniformly. Equation \eqref{eq: lyapunov_1d} is called the Feynman--Kac equation. 
The existence of its solution has been investigated in the literature under different assumptions on $c$. When $c < 0$, the existence of a solution is discussed within the theory of elliptic equations (see, for example, \cite{gilbarg1977elliptic} for details).
Otherwise, one needs to check the principal eigenvalue of the Feynman--Kac operator $\mathcal{L}^Y+\alpha c$ (e.g., using its stochastic representation \cite{liu2002representation,nirenberg}). As shown later in \Cref{sec:construction}, the existence of a solution follows if $Y_t$ satisfies some strong mixing conditions and $\pi^Y(c) < 0$, where $\pi^Y$ is the invariant probability measure of $Y_t$ and $\pi^Y(c)$ is the expectation of $c(y)$ with respect to $\pi^Y$.

The approach described above to construct Lyapunov functions is referred to as \emph{the Feynman--Kac--Lyapunov method}. We extend the framework to a broader class of higher-dimensional SDEs that generalize the driving-driven system in \eqref{toyexample}.
The two key assumptions are that the driven variable $X_t$ has a natural Lyapunov function with a $y$-dependent prefactor and that the driving variable $Y_t$ satisfies SDEs independent of $X_t$ as $X_t$ approaches 0 or infinity. If the system is strictly one-sided coupled as in equation \eqref{toyexample} or the driving system $Y_t$ has a compact phase space, the existence of a Lyapunov function can be rigorously proved. Otherwise, since the perturbation term in the driving dynamics may interact with the Feynman--Kac eigenfunction in a nontrivial way as $|Y_t|$ approaches infinity, one needs to solve the Feynman--Kac equation in a bounded domain and construct a $C^2$ extension outside of the domain. %The existence of a positive solution to the Feynman--Kac equation can be rigorously justified, while the $C^2$ extension usually has to be done numerically. 
We provide statements to construct the Lyapunov function using the Feynman--Kac--Lyapunov method in \Cref{sec:construction}.

As stated above, the Feynman--Kac--Lyapunov method requires a $C^2$ extension of the solution to the Feynman--Kac equation on a bounded domain, which is often completed numerically. In addition, as required by the main theorem proof, the dynamics of the driving variable $Y_t$ must converge to the equation of the $Y$-component as $|X|$ approaches its limit; this requirement can be relaxed if we have an estimate of the solution to the Feynman--Kac equation. Therefore, it is desirable to solve the Feynman--Kac equation numerically for estimation purposes.
In this paper, we develop a numerical method, inspired by \cite{zhai2022deep}, that uses artificial neural networks to solve the Feynman--Kac equation. Since the Feynman--Kac equation has a stochastic representation, we compute approximate solutions at a small number of collocation points using Monte Carlo simulations. We adopt the framework of physics-informed neural networks (PINNs), which represent the solution with neural networks and train them to satisfy the Feynman--Kac equation, as well as the data constraint from the Monte Carlo simulations. This approach does not rely on spatial discretization and therefore applies to high-dimensional problems.

This paper focuses on the construction of the Lyapunov function. The verification of the minorization condition is not our primary focus and is generally carried out through routine calculations. 
All results in this paper, including the corollaries in the examples, can be extended to standard stochastic stability arguments, with additional conditions to provide the minorization condition. For example, if the set $C= \{ y: V (y) \leq c\}$ is ``accessible'' for a constant $c$, which means there exists a point $y \in C$ satisfying $P^{t_0}(x, B_\delta(y)) > 0$ for all $x\in C$, then the minorization condition follows from the hypoelliptic condition, which can be verified by calculating the Lie brackets related to the infinitesimal generator of equation \eqref{SDE1}. 

In \Cref{sec: preliminary}, we outline the problem setting and provide sufficient probabilistic background for Lyapunov functions, stochastic stability, the Feynman--Kac theory, and the multiplicative ergodic theorem. The main results of constructing Lyapunov functions are given in \Cref{sec:construction}. In \Cref{sec:numerics}, we provide a machine learning approach to compute the prefactor $U(y)$ in the Lyapunov function numerically. In \Cref{sec:example}, we showcase several examples of SDEs to construct the Lyapunov functions through analytical or numerical approaches.

\section{Probability preliminaries}\label{sec: preliminary}
Let $\{X_t\}_{t\geq 0}$ be the stochastic process, governed by the SDE 
\begin{equation}\label{eq:SDE_basic}
    \mathrm{d}X_t = b(X_t)\,\mathrm{d}t + \sigma(X_t)\,\mathrm{d}B_t,
\end{equation}
where $X_t,b(X_t)\in \mathbb{R}^{d}$, $\sigma(X_t)\in \mathbb{R}^{d\times n}$, and $B_t$ is the $n$-dimensional Wiener process. If the limit exists, the infinitesimal generator $\mathcal{L}$ of $X_t$ is defined by
\begin{equation}
    \mathcal{L}f(x) = \lim_{t\to 0}\frac{\mathbb{E}_x[f(X_t)] - f(x)}{t},
\end{equation}
where $f\in C_0$ and $\mathbb{E}_x[f(X_t)]=\mathbb{E}(f(X_t)|X_0=x)$. We use the subscript notation with respect to $x$ to indicate the stochastic process starts at $x$, and $\mathcal{L}f(x)$ gives the time derivative of the expectation $\mathbb{E}_x[f(X_t)]$. If $f\in C_0^2$,  then the generator $\mathcal{L}$ is 
\begin{equation}
    \mathcal{L} = \sum_{i=1}^{d}b_i\frac{\partial}{\partial x_i} + \sum_{i,j=1}^{d}a_{ij}\frac{\partial^2}{\partial x_i\partial x_j}
\end{equation}
with $a=\frac{1}{2}\sigma\sigma^T$ being the diffusion matrix.
Throughout this paper, we assume that $\mathcal{L}^Y$ is uniformly elliptic,  meaning the diffusion matrix $a$ is strictly positive definite and has a uniform lower bound, unless specified otherwise.

\subsection{Lyapunov function and stochastic stability} 
For the SDE \eqref{eq:SDE_basic}, a {\it Lyapunov function} is a twice continuously differentiable function $V:\mathbb{R}^{d}\to [1, \infty)$ that satisfies the {\it drift condition}: there exist constants $a>0$ and $b<\infty$ such that
\begin{equation}
    \mathcal{L} V \leq -aV + b.
\end{equation}
The drift condition ensures that the expectation of $V$ has at most linear growth, which in turn provides a means to prove the stability of the stochastic process without requiring an explicit solution of the SDE. 
 
\subsection{Minorization and stochastic stability}
Consider a Markov process $X_t$, a state space $E$, and its Borel $\sigma$-field $\mathcal{B}(E)$. Let $P^t(x, \cdot)$ be the transition kernel of $X_t$, defined as $P^t(x, A) = \mathbb{P}[X_t \in A \,| X_0 = x]$. The {\it minorization condition} states that there exists a constant $\delta>0$, a set $C\in\mathcal{B}(E)$, and a probability measure $\nu$ with $\nu(C)=1$, such that for any $t>0$,
\begin{equation}
P^{t_0}(x,\cdot) \geq \delta \nu(\cdot), \quad \forall x\in C. 
\end{equation}
The minorization condition provides a uniform lower bound on the transition probability on set $C$. As a consequence, the minorization condition enables coupling of the Markov process. 
More precisely, under the minorization condition, a coupling between $(X^{(1)}_t, X^{(2)}_t)$ on $E \times E$ can be constructed in a way such that the two margins are the original Markov process $X_t$ that starts from different initial conditions. Once $X^{(1)}_t$ and $X^{(2)}_t$ enter $C$ simultaneously, the minorization condition allows $\{X^{(1)}_t = X^{(2)}_t\}$ to occur with strictly positive probability, while maintaining two marginal distributions unchanged. As a consequence, the drift condition and the minorization condition imply the convergence to a unique invariant distribution.
We refer to \cite{meyn1993MC, meyn1994convergenceRates, rosenthal1995minorization} for additional details.

\subsection{Feynman--Kac equation and its stochastic representation} 
Let $u(x)$ be the positive solution of the Feynman--Kac equation on a regular domain $\Omega$, which satisfies
\begin{equation} \label{eq:CD}
\begin{aligned}
    &\mathcal{L} u + c_0u = -q, && \quad \text{in } \Omega,\\
        &u = \varphi, && \quad \text{on }\partial \Omega,
\end{aligned}
\end{equation}
for some functions $q,\varphi>0$. The solution to the Feynman--Kac equation \eqref{eq:CD} on a regular domain has a stochastic representation. A point $x\in \partial \Omega$ is said to be {\it regular} if the first exit time of $X_t$ starting from $x$ to leave $\Omega$ is 0 \cite{Oksendal}. This means that $X_t$ starting from $x\in\partial\Omega$ leaves $\Omega$ immediately. 
Let $\Omega\subset \mathbb{R}^d$ be a regular bounded domain, which means that every point $x\in \partial \Omega$ is regular.
Suppose $X_0=x\in \Omega$ and $\tau$ is the first exit time from $\Omega$. Assume
\begin{equation}
    \mathbb{E}_x\left[ \int_{0}^{\tau}\exp\left( \int_{0}^{t}c_0(X_s)\,\mathrm{d}s \right)\,\mathrm{d}t \right] < \infty.
\end{equation}  
If a bounded $C^2$ solution of equation \eqref{eq:CD} exists, then the solution has the stochastic representation \cite{freidlin1985functional}:
\begin{equation}\label{eq:sto_rep}
    u(x) =  \mathbb{E}_x\left[ \int_{0}^{\tau}q(X_t)\exp\left( \int_{0}^{t}c_0(X_s)\,\mathrm{d}s \right)\,\mathrm{d}t \right] + \mathbb{E}_x \left[ \varphi(X_{\tau})\exp\left( \int_{0}^{\tau}c_0(X_s)\,\mathrm{d}s \right) \right].
\end{equation}

The principal eigenvalue of $-(\mathcal{L} + c_0)$ on $\Omega$, denoted by $\lambda_1$, also has a stochastic representation (see Theorem 6.2 of \cite{liu2002representation}) if $\Omega$ has a regular boundary, which is
\begin{equation}
    \lambda_1 = -\lim_{t\to\infty}\frac{1}{t}\ln{\left(\sup_{x\in\Omega}\mathbb{E}_x\left[\mathbf{1}_{\{t<\tau\}}\exp{\left(\int_0^{t}c_0(X_s)\,\mathrm{d}s\right)}\right]\right)}.
\end{equation}
If the principal eigenvalue $\lambda_1 > 0$, the following existence result holds.
\begin{theorem}
  \label{nirenberg}
  Assume $\lambda_{1} >0$ and the boundary condition $\varphi(x) = 0$. Then given $q \in L^{d}(\Omega)$, equation \eqref{eq:CD} has a unique solution $u \in W^{2, d}_{loc}$. In addition, $|u|$ is bounded.
\end{theorem}

\subsection{Multiplicative ergodic theorem}
Consider a generic continuous-time Markov process $X_t$ on $\mathbb{R}^d$. 
\begin{assumption} \label{assume:ab}
    Suppose that 
\begin{enumerate}
    \item[(A)] $X_t$ is geometrically ergodic with an invariant probability measure $\pi$ such that $\pi(V^2) < \infty$ for a Lyapunov function $V$, and 
    \item[(B)] there is a measurable function $F: \mathbb{R}^d\rightarrow [ -1, 1]$ such that $\pi(F) = 0$ and $F$ has nontrivial asymptotic variance with respect to $\pi$.
\end{enumerate}
\end{assumption}
Under \Cref{assume:ab}, Theorem 4.16 in \cite{meyn2003spectral} is known as the multiplicative ergodic theorem. 

\begin{theorem}[Theorem 4.16 in \cite{meyn2003spectral}]
 \label{MET}
Suppose \Cref{assume:ab} holds for $X_t$ and $F$. Then for any complex number $\alpha$ with sufficiently small $|\alpha|$, the operator $\mathcal{L} + \alpha F$ has a maximal and isolated eigenvalue $\Lambda(\alpha)$ with an eigenfunction $\check{f}_{\alpha}(x)$, which satisfy
\begin{equation}
\label{eq:2-10}
\left| \mathbb{E}_{x}\left[ \exp \left ( \alpha \int_{0}^{t} F( X_{s})
    \mathrm{d}s  - t \Lambda (\alpha) \right )\right] - \check{f}_{\alpha}(x)
  \right| \leq B_{0} |\alpha| e^{-b_{0}t} V(x)
\end{equation}
for positive constants $b_{0}$ and $B_{0}$. 
\end{theorem}

\subsection{Eigenvalue perturbation theory} 
If $\mathcal{L}$ is the infinitesimal generator of a stochastic differential equation that is geometrically ergodic, then $\Lambda(\alpha)$ is of order $o(\alpha)$ by the classical eigenvalue perturbation theory. In particular, we have the following lemma.
\begin{lem}
\label{eigperturbation}
    Suppose \Cref{assume:ab} holds for $X_t$ and $F$, then $\lim_{\alpha \rightarrow 0} \Lambda(\alpha)/\alpha = 0$.
\end{lem}
The proof of \Cref{eigperturbation} follows from classical eigenvalue perturbation theory; see \cite{kato2013perturbation}.

\medskip
Let $T$ be a closed linear operator in space $\mathbb{R}^d$. Let $R(\zeta) = (\zeta I - T)^{-1}$ be the resolvent of $T$. The {\it eigenprojection} of $T$ with respect to a semisimple isolated eigenvalue $\lambda$ is denoted by 
\begin{equation}
    \label{eigenproj}
    P= -\frac{1}{2\pi i} \oint R(\zeta) \mathrm{d}\zeta ,
  \end{equation}
where the contour integral integrates $R$ around $\lambda$.

Let $T(\kappa) = T + \kappa T_{1}$ be a small perturbation of $T$. Let $R_{\kappa}$ and $P_{\kappa}$ be the associated resolvent and the eigenprojection with respect to $T(\kappa)$ and $\lambda(\kappa)$. 
The semisimple isolated eigenvalue $\lambda(\kappa)$ is {\it stable} under the perturbation $T(\kappa)$ if (i) the region of convergence, in which $R_{\kappa}(\zeta) \rightarrow R(\zeta)$ as $\kappa \rightarrow 0$ in the generalized sense (see Section 4.2 of \cite{kato2013perturbation}), contains a neighborhood of $\lambda$ with the exception of $\lambda$ itself, and (ii) $\mathrm{Dim} \,P_{\kappa} \leq \mathrm{Dim} \,P$
  for all sufficiently small $\kappa$, where $\mathrm{Dim} P$ means the dimension of the range of $P$. 

Using the notation above, we have the following theorem regarding the eigenvalues of $T(\kappa)$.

  \begin{theorem}[Theorem 8.2.6 of \cite{kato2013perturbation}]
  \label{kato}
  Let $\lambda$ be an isolated, semisimple eigenvalue of a linear operator $T$ with an eigenprojection $P$ such that $\mbox{Dim}\, P = m <\infty$.
  Denote by $D(T_1)$ the domain of $T_1$. Let $\lambda$ be stable and let $P \mathbb{R}^d \subset D(T_{1})$. Then the eigenvalues of $T(\kappa)$ admit asymptotic expansions
$$
  \mu_{j}(\kappa) = \lambda + \kappa \mu_{j}^{(1)} + o( \kappa) ,
  \quad j = 1, \ldots, m ,
$$
where $\mu_{j}^{(1)}$ are the repeated eigenvalues of the operator $PT_{1}P$ considered in the $m$-dimensional space $PX$.
\end{theorem}

\begin{proof}[Proof of \Cref{eigperturbation}] 
Because $F$ and $X_t$ satisfy \Cref{assume:ab}, \Cref{MET} shows that the maximal eigenvalue of $\mathcal{L} + \alpha F$, denoted by $\Lambda( \alpha)$, is isolated for all sufficiently small $|\alpha|$. Therefore, $0$ must be an isolated maximal eigenvalue of $\mathcal{L}$. 

Since $F$ is bounded, $T_{1}$ is a bounded linear operator. Therefore, by Theorem 4.3.11 of \cite{kato2013perturbation}, the resolvent $R_{\alpha}$ is continuous under small perturbations. By Theorem 4.3.16 of \cite{kato2013perturbation}, the spectrum is upper semicontinuous with respect to small perturbations (meaning that $\mathrm{Dim} P_{\alpha} = \mathrm{Dim} P$ for all small $\alpha$). Hence $0$ is a stable eigenvalue by the definition of stable eigenvalue given above. 

Let $P$ be the eigenprojection of $T$ with
respect to eigenvalue $0$. It remains to study the properties of $P$. Since $\pi$ is unique and geometrically ergodic, the multiplicity of the eigenvalue $0$ must be $1$. Thus, $\mathrm{Dim}  P=1$. Therefore, the eigenfunction of $T$ with respect to $0$ is a constant function $\mathbf{1}$. The eigenfunction of $T^*$ with respect to $0$ is the probability density function of $\pi$, denoted by $u$. In addition, denoting the resolvent of $T^*$ by $R^*$, it follows that $R^*(\bar{\zeta}) = (R(\zeta))^*$. Therefore, since $0$ is a real number, the eigenprojection of $T^*$ with respect to $0$ must be $P^*$. 

Note that $P$ is the projection into the set of constant functions; for any function $w$ we have $P w = C_w \mathbf{1}$. Integrating $Pw$ with $u$, we have
$$
\langle u, Pw \rangle = C_w \langle u, \mathbf{1} \rangle = C_w
$$
because $\langle u, \mathbf{1} \rangle = 1$. Therefore, we have
$$
C_w = \frac{\langle u, Pw\rangle}{\langle u, \mathbf{1} \rangle} = \frac{\langle  P^* u ,w\rangle}{\langle u, \mathbf{1}\rangle} = \langle u, w\rangle
$$
because $P^*$ projects $u$ into itself. This gives $P w = \pi (w)\mathbf{1}$. Therefore, we have $P T_1 P w = \pi(w)P T_1 \mathbf{1} = \pi(w)PF = \pi(w)\pi(F) = 0$ because $T_1 \mathbf{1} = F$ and $\pi(F) = 0$. Hence, by Theorem \ref{kato}, we have
$$
  \Lambda(\alpha) = o(\alpha) .
  $$    
\end{proof}

\subsection{Existence of a solution} 
Theorem 2.7 in \cite{hani2025nonequilibrium} guarantees that there exists a solution of the Feynman--Kac equation  \eqref{eq:CD} if $c_0$ is a sufficiently small function. More precisely, we consider a generic $C^3$ function $c(x)$ and let $c_0 = \alpha c$. The goal is to show that equation \eqref{eq:CD} has a strictly positive solution when $\alpha$ is sufficiently small. We introduce the proof for the sake of completeness of this paper. 

If $c$ is defined globally, we let $\hat{c}= \beta\mathbf{1}_{\Omega}c$ and shift $\hat{c}$ by its mean as $F := \hat{c}-\pi(\hat{c})$. For any function $c\in C^3$, we choose the constant $\beta>0$ such that the range of $F$ is contained in $[-1,1]$.

\begin{theorem}[Theorem 2.7 of \cite{hani2025nonequilibrium}]
    \label{NLS27}
    Let $\Omega\subset\mathbb{R}^d$ be a regular bounded domain. Let $c$ be a $C^3$ function, and $F$ be the rescaled function of $c$ as above. If \Cref{assume:ab} holds for $X_t$ and the observable $F$, and $\pi(\hat{c}) < 0$, then the following Feynman--Kac equation has a positive solution for all $\alpha \in (0, \alpha^*)$,
\begin{equation}\label{eq:NLS27}
    \begin{aligned}
        &\mathcal{L}u + \alpha cu = -q, && \quad  \text{in }\Omega,\\
        &u = \varphi, && \quad \text{on }\partial \Omega,
    \end{aligned}
\end{equation}
where $\alpha^* > 0$, $q>0$ is a bounded $C^2$ function, and $\varphi\in C^2$ is positive and bounded. 
\end{theorem}
\begin{proof}
Recall that $F = \hat{c}-\pi(\hat{c})$ and $\pi(\hat{c}) < 0$. By \Cref{MET}, we have
\begin{equation}
    \mathbb{E}_{x}\left[\exp\left(\alpha_0\int_0^t \hat{c}\left(X_s\right)\,\mathrm{d}s\right)\right] \leq \left[\check{f}_{\alpha_0}(x)+B_0|\alpha_0|e^{-b_0t}V(x)\right]\exp\left(t\Lambda(\alpha_0) + t\alpha_0\pi(\hat{c})\right).
\end{equation}
In addition, the eigenvalue perturbation result in \Cref{eigperturbation} implies that $\lim_{\alpha_0\to 0}\Lambda(\alpha_0)/\alpha_0=0$. Note that the eigenvalue perturbation result only depends on the operator $\mathcal{L}$ defined in $\mathbb{R}^d$ and $F$. Therefore, for all sufficiently small $\alpha_0$, we have
\begin{equation}
\mathbb{E}_{x}\left[\exp\left(\alpha_0\int_0^t \hat{c}\left(X_s\right)\,\mathrm{d}s\right)\right] \leq B_1\exp\left(-\alpha_0 b_1 t\right)
\end{equation}
for some constant $B_1<\infty$ and $b_1>0$. The right-hand-side bound is uniform for any $x \in \Omega$.  By the stochastic representation of the principal eigenvalue of the Feynman--Kac operator (Theorem 6.2 of \cite{liu2002representation}), this implies that the principal eigenvalue of $-(\mathcal{L} + \alpha_0 \hat{c})$  in equation \eqref{eq:NLS27} satisfies $0 < \alpha b_1 \le \lambda_1$. Thus, letting $\alpha = \alpha_0 \beta$, equation \eqref{eq:NLS27} with $\varphi=0$ has a solution in $W^{2, d}_{loc}$ by Theorem 1.2 of \cite{nirenberg}. Let $w$ be a bounded $C^2$ function in $\Omega$ such that $w=\varphi$ on $\partial \Omega$. By letting $\hat{u}=u-w$, we can rewrite \eqref{eq:NLS27} as a problem with zero boundary. Therefore, equation \eqref{eq:NLS27} has a bounded solution in $W^{2, d}_{loc}$. 

By using the higher regularity of solutions to elliptic equations (see, for example, Theorem 6.5 of \cite{evans2022partial}), the solution satisfies $u \in C^2(\Omega)$ if all coefficients possess sufficient regularity. In this case, the positivity of $u$ follows directly from the stochastic representation in equation \eqref{eq:sto_rep}. Therefore, $u$ must be strictly positive.
\end{proof}

\section{Construction of Lyapunov functions}\label{sec:construction}
In this section, we provide the main theorem that describes how to construct a Lyapunov function for the driving-driven SDEs. 
\subsection{Main theorem}
Consider the following SDE for a driving-driven pair $(X_t, Y_t)$ with $X_t\in\mathbb{R}^{d_1}$, $Y_t\in\mathbb{R}^{d_2}$:
\begin{equation}
\label{SDE2}
    \begin{aligned}
        \mathrm{d}X_t &= f(X_t, Y_t)\,\mathrm{d}t + \sigma_x(X_t) \mathrm{d}B^x_t,\\
        \mathrm{d}Y_t &= g(X_t, Y_t)\,\mathrm{d}t + \sigma_y(Y_t)\,\mathrm{d}B^y_t,
    \end{aligned}
\end{equation}
where $\sigma_x$ and $\sigma_y$ are matrix-valued functions in $\mathbb{R}^{d_1\times n_1}$ and $\mathbb{R}^{d_2\times n_2}$, respectively, and $B^x_t\in\mathbb{R}^{n_1}$ and $B^y_t\in\mathbb{R}^{n_2}$ are the Wiener processes. Let $\mathcal{L}$ be the infinitesimal generator for $(X_t,Y_t)$ for equation \eqref{SDE2}. 

Throughout this paper, we assume all coefficient functions in $\mathcal{L}$ are in $C^3$. 
In addition, we assume that the vector field $g(x,y)$ can be approximated by a vector field $g_0(y)$ that only depends on $y$. In addition, their difference is bounded by
\begin{equation}
    \label{approximation}
    |g(x,y)-g_0(y)|\leq h(x,y),
\end{equation}
for all $x,y$, where $|\cdot|$ denotes the 2-vector norm throughout this paper.
We define the infinitesimal generators $\mathcal{L}^{X}$ and $\mathcal{L}^{Y_0}$ with respect to the vector fields $f$ and $g_0$, respectively, which are
\begin{equation}
\begin{aligned}\mathcal{L}^{X}w(x) &= \nabla w(x)\cdot f(x,y) +\frac{1}{2}\sum_{i,j=1}^{d_1}\left(\sigma_x\sigma_x^T\right)_{ij}\frac{\partial^2 w(x)}{\partial x_i\partial x_j}, \\
\mathcal{L}^{Y_0}v(y) &= \nabla v\cdot g_0(y)+\frac{1}{2}\sum_{i,j=1}^{d_2}\left(\sigma_y\sigma_y^T\right)_{ij}\frac{\partial^2 v(y)}{\partial y_i\partial y_j},
\end{aligned}
\end{equation}
for test functions $w\in C^2(\mathbb{R}^{d_1})$ and $v\in C^2(\mathbb{R}^{d_2})$. 
To prepare for the main theorem, we make the following assumptions about equation \eqref{SDE2}.

\begin{assumption}\label{general_assumption}
Equation \eqref{SDE2} satisfies the following properties:
\begin{enumerate}
    \item[(A)] There exists a family of Lyapunov functions $V_\alpha(x) \geq 1$ and a function $c(y)$, such that for sufficiently small $\alpha > 0$, the inequality
    \begin{equation} \label{eq:assumeA}
        \mathcal{L}^X V_\alpha(x)\leq\alpha c(y)V_\alpha(x)
    \end{equation}
    holds for all $x$ with sufficiently large $V_\alpha$-value. 
    \item[(B)] There exists a pair of bounded regular domains $\Omega_\text{I}$ and $\Omega_\text{O}$, with
      $\Omega_\text{I} \subset \Omega_\text{O}\subseteq\mathbb{R}^{d_2}$, an exterior Lyapunov function
      $\tilde{U} \geq 1$, a nonnegative $C^2$ function $q(y) \geq 0$, a uniformly positive function $p(y)$, and a constant $\zeta > 0$, such that 
    \begin{equation}\label{eq: assumption_U_in}
    \begin{aligned}
          & \mathcal{L}^{Y_0}\tilde{U}(y)\leq-p(y)\tilde{U}(y), && \forall y\notin\Omega_\text{I} \\ 
          & \mathcal{L}^{Y_0}\tilde{U}(y) + \alpha c(y) \tilde{U}(y) \leq q(y),  &&\forall y \in \Omega_\text{I} \\
   & p(y) > \zeta c(y), && \forall y \not\in \Omega_\text{I} \\
   & p(y) > \zeta \frac{|\nabla \tilde{U}|(y)}{\tilde{U}(y)},  && \forall y \not\in \Omega_\text{I}
   \end{aligned}
\end{equation}
    for all sufficiently small $\alpha$ and $c(y)$ in Assumption (A).
    \item[(C)] Let $\varphi \in C^3$ be a positive function defined on $\partial \Omega_\text{I}$. There exists a constant $  \gamma> 0$ such that the Feynman--Kac equation
    \begin{equation}
    \label{eq:feynman}
        \begin{aligned}
            &\mathcal{L}^{Y_0}Q+\alpha cQ = -q(y) - \gamma, && \quad \text{ in }\Omega_\text{I}, \\
            & Q = \varphi, && \quad \text{ on }\partial\Omega_\text{I},
        \end{aligned}
    \end{equation}
    has a positive solution $Q$ in $C^2(\Omega_\text{I})$ for $\alpha$, $c(y)$, and $q(y)$ in Assumptions (A) and (B).
    \item[(D)] For some $\alpha\in (0,\zeta/3)$ and a positive constant $C$, the solution $Q(y)$ to the Feynman--Kac equation \eqref{eq:feynman} has a $C^{2}$ extension of $Q$ in Assumption (C) to $\mathbb{R}^{d_{2}}$ such that
\begin{equation}
  \label{eq:extension}
\begin{aligned}
  &\mathcal{L}^{Y_{0}} Q + \alpha c Q \leq \frac{1}{2} p \tilde{U}, && y \in \Omega_\text{O} \setminus \Omega_\text{I},\\
  &Q = C, && y \notin \Omega_\text{O}, \\
  &\tilde{U} \geq C, && y \in \Omega_\text{O}^c.
\end{aligned}
\end{equation}
    \item[(E)] Let $h(x,y)$ be the difference given in \eqref{approximation}, then $\lim\limits_{|x|\to\infty}h(x,y)=0$ uniformly for all $y \in \Omega_\text{I}$. 
\end{enumerate}
\end{assumption}

\Cref{general_assumption}(C) requires that the Feynman--Kac equation has a strictly positive $C^2$ solution. In this paper, we use both analytical and numerical approaches to verify the existence of a positive solution for the Feynman--Kac equation. 
From an analytical perspective, notice that the driving variable $Y_t$ is approximated by $\tilde{Y}_t$, which satisfies
    \begin{equation}
        \label{eq:y0}
        \mathrm{d}\tilde{Y}_t = g_0(\tilde{Y}_t) \,\mathrm{d}t + \sigma_y(\tilde{Y}_t) \,\mathrm{d}B^y_t,
    \end{equation}
when $h$ is sufficiently small by assumption. Let $\pi^{\tilde{Y}}$ be the unique invariant probability measure of $\tilde{Y}_t$ in \eqref{eq:y0}. If the SDE \eqref{eq:y0} and $\pi^{\tilde{Y}}$ satisfy conditions of \Cref{NLS27}, then
    \Cref{general_assumption}(C) holds for sufficiently small $\alpha
    > 0$.
    
From a numerical perspective, we solve the Feynman--Kac equation \eqref{eq:feynman} numerically and examine whether the solution is strictly positive to verify \Cref{general_assumption}(C). After obtaining a numerical solution $Q$ to the Feynman--Kac equation within the domain $\Omega_\text{I}$, we numerically extend $Q$ from $\Omega_\text{I}$ to the domain $\Omega_\text{O}$ to verify \Cref{general_assumption}(D). 
In addition, one may further relax \Cref{general_assumption}(E), which we illustrate in \Cref{sec:numerics}.

\begin{theorem}\label{thm:ly}
    Suppose that \Cref{general_assumption} holds for the driving-driven system \eqref{SDE2}, then there exists a strictly positive Lyapunov function $W(x,y)$ and constants $\alpha,\delta>0$ such that for every $x$ with sufficiently large $V_{\alpha}(x)$, we have $\mathcal{L}W(x,y)\leq -\delta W(x,y)$ for all $y\in\mathbb{R}^{d_2}$.
\end{theorem}
\begin{proof}
    Define $W(x,y)=V_{\alpha}(x)U(y)$, where $V_{\alpha}$ is the $x$-Lyapunov function
    in \Cref{general_assumption} and $U = \tilde{U} + Q$ is a prefactor function, where $\tilde{U}$ and $Q$ are from \Cref{general_assumption}(B) and \Cref{general_assumption}(C). 
    Direct computation yields
    \begin{equation}\label{LW}
   \mathcal{L}W(x,y)=U(y)\mathcal{L}^X V_{\alpha}(x)+V_{\alpha}(x)\mathcal{L}^YU(y),
    \end{equation}
    where $\mathcal{L}^X $ and $\mathcal{L}^Y$ are the $x$-component and $y$-component of $\mathcal{L}$, respectively. 

      \medskip

      {\bf Case 1}: $y \in \Omega_\text{I}$. If $y \in \Omega_\text{I}$, we consider the Feynman--Kac equation for the Lyapunov function of the $y$ variable, which is equation \eqref{eq:feynman}. 
    By \Cref{general_assumption}(C), equation \eqref{eq:feynman} has a strictly positive $C^2$ solution $Q(y)$ for all sufficiently small
    $\alpha > 0$. 

      By \Cref{general_assumption}(E), we have $g(x,y)$ as a perturbation of $g_0(y)$ with the difference bound $h(x,y)$. Therefore, direct computations yield 
    \begin{equation}
    \begin{aligned}
        \mathcal{L}^YU&=\nabla U\cdot g(x,y)+\frac{1}{2}\sum_{i,j=1}^{d_2}\left(\sigma_y\sigma_y^T\right)_{ij}\frac{\partial^2 U}{\partial y_i \partial y_j}\\
        &\leq \nabla U\cdot g_0+\frac{1}{2}\sum_{i,j=1}^{d_2}\left(\sigma_y\sigma_y^T\right)_{ij}\frac{\partial^2 U}{\partial y_i \partial y_j}+h(x,y)|\nabla U|\\
        &=\mathcal{L}^{Y_0}U+h(x,y)|\nabla U|,
      \end{aligned}
    \end{equation}
where $|\nabla U|(y)$ is the 2-vector norm of $\nabla U$ at $y$.
    
Using $U = \tilde{U} + Q$, we compute $\mathcal{L}W$ for $y\in\Omega_\text{I}$
    \begin{equation} \label{eq: LWxy2}
    \begin{aligned}
        \mathcal{L}W(x,y) \leq & %U\alpha c(y)V+V\mathcal{L}^{Y_0}U+Vh(x,y)|\nabla U|\\
        \left[\mathcal{L}^{Y_0}U(y)+\alpha c(y)U(y)+h(x,y)|\nabla
      U(y)|\right]V_{\alpha}(x)\\
      \leq & \left[\mathcal{L}^{Y_0}\tilde{U}(y)+\alpha
             c(y)\tilde{U}(y)+\mathcal{L}^{Y_0}Q(y)+\alpha
             c(y)Q(y)\right]V_{\alpha}(x) + h(x,y) |\nabla
             U| V_{\alpha} \\
      \leq & (q -q - \gamma )V_\alpha +  h(x,y) |\nabla
             U| V_{\alpha},
      \end{aligned}
      \end{equation}
      where we use the assumptions in \eqref{eq:assumeA}, \eqref{eq: assumption_U_in}, and \eqref{eq:feynman}.
      Since $\Omega_\text{I}$ is bounded and function 
      $|\nabla U |$ is continuous, by \Cref{general_assumption}(E), there exists a sufficiently large $|x|$ such that   
\begin{equation}  
    h(x,y) |\nabla U|(y) < \frac{1}{2} \gamma
\end{equation}
  uniformly for all $y \in \Omega_\text{I}$.  This implies
\begin{equation}
  \mathcal{L}W(x,y) \leq  -\frac{1}{2} \gamma V_{\alpha}  < -\delta W
\end{equation}
  for some constant $\delta > 0$, because $U = \tilde{U} + Q$ is bounded in $\Omega_\text{I}$.

  \medskip
 
  {\bf Case 2}: $y \notin \Omega_\text{I}$. If $y\notin\Omega_\text{I}$, by
  \Cref{general_assumption}(D), the solution $Q(y)$ to the Feynman--Kac equation has a global $C^{2}$ extension for some $\zeta/3 > \alpha > 0$, such that $\mathcal{L}^{Y_{0}} Q + \alpha c Q \leq \frac{1}{2} p(y) \tilde{U}$ for all $y\in \Omega_\text{O}\setminus\Omega_\text{I}$. In addition, we have $\mathcal{L}^{Y_{0}}Q = 0$ and $\mathcal{L}^{Y_{0}}\tilde{U} \leq - p(y) \tilde{U}$ if $y\notin \Omega_\text{O}$. This gives
    \begin{equation}
    \begin{aligned}
        \mathcal{L}W(x,y) &\leq \left ( \mathcal{L}^{Y_{0}}U(y) + \alpha c(y)U(y) +h(x,y) |\nabla U|(y) \right ) V_{\alpha}(x) \\
        \leq & \left (-\frac{1}{2} p(y) \tilde{U} + \alpha c(y)\tilde{U}(y) + h(x,y)|\nabla U|(y) \right ) V_{\alpha}(x) .
    \end{aligned}
  \end{equation}
If $y \in \Omega_\text{O}\setminus \Omega_\text{I}$, because $\alpha c(y) \leq p(y)/3$ and $h(x,y)|\nabla U|(y)$ are uniformly bounded, we can make $|x|$ sufficiently large such that
\begin{equation}
  \alpha  c(y)U(y) + h(x,y)|\nabla U|(y) < \frac{1}{8} p(y) \tilde{U} .
  \end{equation}
  This gives
  \begin{equation}
   \mathcal{L}W(x,y)  \leq -\frac{1}{8} p(y)
  \tilde{U} V_{\alpha} \leq - \delta W 
  \end{equation}
 for some constant $\delta > 0$. This is because $\Omega_\text{O}$ is bounded, so there must be $\delta > 0$ making $\frac{1}{8} p(y)\tilde{U} \geq \delta U$. 

 If instead $y \notin \Omega_\text{O}$, since $Q = C$ is a constant, now we have
 \begin{equation}
\begin{aligned}
  \mathcal{L}W(x,y) &\leq \left (-p(y) \tilde{U} + \alpha
        c(y)(\tilde{U}(y) + C) + h(x,y)|\nabla \tilde{U}|(y) \right
                      ) V_{\alpha}(x) \\
  &\leq V_{\alpha} \tilde{U} \left( -p(y) + 2\alpha c(y)+ h(x,y) \frac{|\nabla \tilde{U} |}{
     \tilde{U}}\right)\\
  &\leq   V_{\alpha} \tilde{U} \left ( -\frac{1}{3}p(y) + h(x,y) \frac{|\nabla \tilde{U} |}{
         \tilde{U}} \right )  \\
    &\leq p(y)V_{\alpha}\tilde{U}\left(-\frac{1}{3}+h(x,y)\zeta^{-1}\right)\\
  &\leq  - \frac{1}{4} p(y) V_{\alpha} \tilde{U}
\end{aligned}
\end{equation}
for sufficiently small $h$. Since $\tilde{U} \geq
U/2$ and $p(y)$ is uniformly positive in $\Omega_\text{O}^{c}$, we have
\begin{equation}
   \mathcal{L}W(x,y) \leq - \frac{1}{4} p(y) V_{\alpha} \tilde{U}  \leq -
   \delta W
   \end{equation}
   for some strictly positive constant $\delta > 0$.

   Combining all cases above, we can choose $\alpha$ sufficiently small such that
\begin{equation}
  \mathcal{L}W \leq - \delta W
\end{equation}
for some $\delta > 0$ when $|x|$ is large enough to make $h(x,y)$ uniformly small. This completes the proof. 
\end{proof}

Notice that \Cref{thm:ly} guarantees the drift condition for the function $W(x,y)$ when $V_\alpha(x)$ is large. Whether $W(x,y)$ is a suitable Lyapunov function that leads to stochastic stability also depends on the $Y$-dynamics and the property of $V_{\alpha}(x)$. For example, if $V_{\alpha}(x) \rightarrow \infty$ as $|x| \rightarrow \infty$, $\mathcal{L}^{Y}\tilde{U}$ gives sufficient negative drift for all moderate $V_\alpha(x)$ and large $y$, and $\{ V_{\alpha}(x) \leq M \}$ satisfies the minorization condition for a sufficiently large $M$, then the stochastic stability follows. 

\subsection{Special case 1: one-way coupled system} In the proof of \Cref{thm:ly}, the assumption about the Feynman--Kac equation on a bounded domain mainly comes from the perturbation term $h(x,y)$. When $h(x,y) = 0$ and system \eqref{SDE2} becomes one-way coupled, the assumption can be dramatically simplified. In fact, we have the following theorem for the following SDE
\begin{equation}
\label{SDE3}
    \begin{aligned}
        \mathrm{d}X_t &= f(X_t, Y_t)\,\mathrm{d}t + \sigma_x(X_t) \mathrm{d}B^x_t\\
        \mathrm{d}Y_t &= g(Y_t)\,\mathrm{d}t + \sigma_y(Y_t)\,\mathrm{d}B^y_t .
    \end{aligned}
\end{equation}
Denote the infinitesimal generator of $Y_t$ by $\mathcal{L}^Y$. We have the following assumptions about equation \eqref{SDE3}.
\begin{assumption}\label{general_oneside}
Equation \eqref{SDE3} satisfies the following properties:
\begin{enumerate}
    \item[(A)] There exists a family of Lyapunov functions $V_\alpha(x) \geq 1$ and a function $c(y)$, such that the inequality
    \begin{equation}
        \mathcal{L}^X V_\alpha(x)=\nabla V_{\alpha}(x)\cdot f(x,y)+\frac{1}{2}\sum_{i,j=1}^{d_1}\left(\sigma_x\sigma_x^T\right)_{ij}(V_\alpha)_{ij}(x)\leq\alpha c(y)V_\alpha(x), \quad \forall x \in \mathbb{R}^{d_1}
    \end{equation}
    holds for all sufficiently small $\alpha > 0$.
    \item[(B)] $Y_t$ is geometrically ergodic with an invariant probability measure $\pi$ such that $\pi( V_y^2) < \infty$ for a Lyapunov function $V_y$.
    \item[(C)] $c(y)$ is uniformly bounded and has asymptotic variance with respect to $\pi$. In addition, $\pi(c) < 0$.
\end{enumerate}
\end{assumption}
\begin{theorem}
    \label{thm:oneside}
    Assume \Cref{general_oneside} holds for the one-way coupled system \eqref{SDE3}, then there exists a strictly positive Lyapunov function $W(x,y)$ and constants $\alpha, \delta > 0$, such that for every $x$ with sufficiently large $V_\alpha(x)$, we have $\mathcal{L}W(x,y) \leq - \delta W(x,y)$ for all $y \in \mathbb{R}^{d_2}$.
\end{theorem}
\begin{proof}
Let $\hat{c} := \beta c$ and $F(y) := \hat{c} - \pi( \hat{c})$. By \Cref{general_oneside}, $c$ is uniformly bounded. Therefore, we can choose a suitable $\beta$ such that the range of $F$ is contained in $[-1 , \,1]$. Therefore, \Cref{assume:ab} is satisfied by $Y_t$ and $F$. By \Cref{MET}, the operator $\mathcal{L}^Y + \alpha F$ has a maximal and isolated eigenvalue $\Lambda(\alpha)$ and an eigenfunction $U_\alpha(y)$. Note that $U_\alpha$ must be strictly positive everywhere by \eqref{eq:2-10}. It remains to show that $U_\alpha$ has $C^2$ regularity. Note that when $V$ is bounded, Theorem \Cref{MET} further implies that $U_\alpha$ must be continuous on compact domain, because 
\[
u_t(x) := \mathbb{E}_x[\exp \left ( \alpha \int_0^t F(X_s) \mathrm{d}s - t \Lambda(\alpha)\right)]
\]
converges to $U_\alpha$ uniformly on compact set as $t\rightarrow \infty$. Then we need the following higher regularity argument to show that $U_\alpha$ is $C^2$.

Let $R_r = \int_0^\infty e^{-rt}P_t \mathrm{d}t$ be the diffusion resolvent. Notice that $U_\alpha$ satisfies $U_\alpha=R_r[(r-\Lambda(\alpha)+\alpha F)U_\alpha] := R_r g$. Now let \(B\) be a ball compactly contained in the state space and \(\tau_B\) its first exit time. Splitting the resolvent at \(\tau_B\) using the strong Markov property gives

\begin{align}
\label{eq:u_alpha}
    U_\alpha(x)&= \mathbb E_x\!\left[e^{-r\tau_B}U_\alpha(X_{\tau_B})\right] +\mathbb E_x\!\int_0^{\tau_B}e^{-rs}g(X_s)\,ds\\
&:=H + G, \qquad x\in B.
\end{align}

Because $U_\alpha \in L^\infty_V$, it must be bounded in $B$. The elliptic Dirichlet theory implies that the equation
\[
(r - \mathcal{L}^{Y_0})H_{PDE} = 0 \quad \mbox{in} \quad B; \quad  H_{PDE} = U_\alpha \quad \mbox{on} \quad \partial B\,.
\]
admits a $C^2$ solution and that equation
\[
(r - \mathcal{L}^{Y_0})G_{PDE} = g\quad \mbox{in} \quad B; \quad  G_{PDE} = 0 \quad \mbox{on} \quad \partial B\,.
\]
admits a solution in $W^{2,p}(B)$, with $p>d$. The Itô–Krylov formula then implies that $H = H_{PDE}$ and $G = G_{PDE}$. Therefore, we have $U_\alpha \in W^{2,p}_{loc}$. Taking \(p>d\), Sobolev embedding makes \(U_\alpha\) locally Hölder, and interior Schauder regularity (see \cite{gilbarg1977elliptic}) yields \(U_\alpha\in C^{2,\beta}_{\mathrm{loc}}\) for some \(\beta>0\).

By the eigenvalue perturbation theorem \cite{kato2013perturbation}, we know that $\Lambda(\alpha)$ is $o(\alpha)$ small. Therefore, letting the operator $\mathcal{L}^Y + \alpha c = \mathcal{L}^Y + \alpha(F/\beta + \pi(c))$ act on $U_{\alpha/\beta}$, we have
$$
(\mathcal{L}^Y + \alpha c)U_{\alpha/\beta} = (\Lambda(\alpha/\beta) + \alpha\pi(c)) U_{\alpha/\beta} < 0
$$
for all sufficiently small $\alpha > 0$.

Defining $W(x,y) = V_\alpha(x)U_{\alpha/\beta}(y)$, we have
\begin{align*}
    \mathcal{L}W &= U_{\alpha/\beta}(y) \mathcal{L}^XV_\alpha + V_\alpha(x) \mathcal{L}^Y U_{\alpha/\beta}(y) \leq  V_\alpha(x) \left [ \mathcal{L}^Y U_{\alpha/\beta}(y) + \alpha c(y) U_{\alpha/\beta}(y)\right ] \\
    &= (\Lambda(\alpha/\beta) + \alpha\pi(c)) V_\alpha(x) U_{\alpha/\beta}(y) = (\Lambda(\alpha/\beta) + \alpha\pi(c)) W(x,y) .
\end{align*}
Let $\alpha$ be sufficiently small such that $\Lambda(\alpha/\beta) + \alpha\pi(c) < 0$. The proof is completed by letting $\delta = -\Lambda(\alpha/\beta) + \alpha\pi(c)$.
\end{proof}

\subsection{Special case 2: Power-law Lyapunov function $V$}
In many situations when we study the stochastic stability or the persistence of a system, the natural Lyapunov function is either $|x|^\alpha$ or $|x|^{-\alpha}$, where $| \cdot|$ denotes the 2-norm. In this subsection, we consider the above two special cases of $V_\alpha(x)$ in the context of the following SDE with $X_t\in\mathbb{R}^{d_1}$ and $Y_t \in \mathbb{R}^{d_2}$:
\begin{equation}\label{eq:sp}
    \begin{aligned}
        \mathrm{d}X_t &= f(X_t, Y_t)\,\mathrm{d}t + \varepsilon X_t \odot\,\mathrm{d}B^x_t\\
        \mathrm{d}Y_t &= g(X_t, Y_t)\,\mathrm{d}t + \sigma(X_t, Y_t)\,\mathrm{d} B^y_t,
    \end{aligned}
\end{equation}
where $\varepsilon$ is a positive constant, $\odot$ is the element-wise product, $\sigma$ is a matrix-valued function in $\mathbb{R}^{d_2\times n_2}$, and $B^x_t\in\mathbb{R}^{d_1}$ and $B^y_t\in\mathbb{R}^{n_2}$ are the Wiener processes.
We provide corresponding modified assumptions to construct the Lyapunov function $W(x,y)$ for $(X_t,Y_t)$ in \eqref{eq:sp} for the following two scenarios: for sufficiently small $|x|$ and for sufficiently large $|x|$. We provide the corresponding assumptions in \Cref{assumption1} and \Cref{assumption2} for these two cases, respectively. 

\begin{assumption}\label{assumption1} Equation \eqref{eq:sp} satisfies the following conditions:
\begin{enumerate} 
    \item[(A)] For sufficiently small $|x|$, $- x\cdot f(x,y) \leq \hat{c}(y) |x|^2$ for some function $\hat{c}(y)$.
    \item[(B)] There exists a pair of bounded regular domains $\Omega_\text{I}\subset\Omega_\text{O}\subseteq\mathbb{R}^{d_2}$ and an exterior Lyapunov function $\tilde{U}\geq 1$ such that 
    \begin{equation}
          \mathcal{L}^{Y_0}\tilde{U}(y)\leq-p(y)\tilde{U}(y), \quad \forall y\notin\Omega_\text{I},
    \end{equation}
    and
    \begin{equation}
        \mathcal{L}^{Y_0}\tilde{U}(y) + \alpha\bar{c}(y)\tilde{U}(y)\leq q(y), \quad \forall y\in\Omega_\text{I}
    \end{equation}
    for all sufficiently small $\alpha>0$, some positive $C^2$ function $q(y)$, and some uniformly positive function $p(y)$ that satisfies
     \begin{equation}
            p(y) > \zeta \bar{c}(y) \quad \mbox{and} \quad p(y) > \zeta \frac{|\nabla \tilde{U}(y)|}{\tilde{U}(y)}
    \end{equation}
    for some constant $\zeta>0$ uniformly for $y\in\Omega_\text{O}$ and $\bar{c}(y) = \hat{c}(y) + \frac{1}{2}\varepsilon^2(\alpha + 1)$.
        
    \item[(C)]  The Feynman--Kac equation
    \begin{equation}
    \label{FKsmall}
    \begin{cases}
        \begin{aligned}
            &\mathcal{L}^{Y_0}Q + \alpha \bar{c} Q = -q(y)-\gamma, && \text{in }\Omega_\text{I}\\
            &Q = \varphi, && \text{on } \partial\Omega_\text{I}
        \end{aligned}
    \end{cases}
        \end{equation}
    has a positive solution in $C^2(\Omega_\text{I})$ for all $\alpha>0$, where $\varphi\in C^3$ is a positive function defined on $\partial\Omega_\text{I}$ and $\gamma>0$ is a constant.
    \item[(D)] For some $\zeta/3>\alpha>0$, the solution $Q(y)$ to equation \eqref{FKsmall} has a $C^2$ extension to the whole domain such that 
    \begin{equation}
    \begin{cases}
        \begin{aligned}
            &\mathcal{L}^{Y_0}Q + \alpha\bar{c} Q\leq \frac{1}{2}p(y)\tilde{U}(y), && y\in\Omega_\text{O}\setminus\Omega_\text{I}\\
            &Q = C, && y\notin\Omega_\text{O}
        \end{aligned}
    \end{cases}
        \end{equation}
    for some constant $C$ such that $\tilde{U}\geq C$ on $\Omega_\text{O}^c$.
    \item[(E)] $\lim_{|x|\to 0}h(x,y)= 0$, where $h(x,y)$ is given in \eqref{approximation}.
\end{enumerate}    
\end{assumption}

\begin{cor}
\label{cor1}
 \Cref{assumption1} implies that $|x|^{-\alpha}U(y)$ is a Lyapunov function for equation \eqref{eq:sp} provided $|x|$ is sufficiently small, where $U = \tilde{U} + Q$. 
\end{cor}
\begin{proof}
    Define $W(x,y) = |x|^{-\alpha} U(y)$, where $U=\tilde{U}+Q$ is a prefactor. Therefore,
    \begin{equation}
        \mathcal{L}W(x,y) = U(y) \mathcal{L}^X |x|^{-\alpha} + |x|^{-\alpha} \mathcal{L}^Y U(y).
    \end{equation}
    Direct calculation yields
    \begin{equation}
    \begin{aligned}
            \mathcal{L}^X|x|^{-\alpha} &= -\alpha |x|^{-\alpha - 2}\left[x \cdot f(x,y)\right] + \frac{1}{2}\varepsilon^2 \alpha |x|^{-\alpha}\left [ (\alpha + 2)\frac{\sum_{i = 1}^{d_1} x_i^4}{|x|^4} - 1\right] \\
            &\leq \alpha |x|^{-\alpha}\hat{c}(y) + \frac{1}{2}\varepsilon^2\alpha(\alpha + 1) |x|^{-\alpha}\\
            &=\alpha |x|^{-\alpha} \bar{c}(y),
            \end{aligned}
    \end{equation}
where $\bar{c}(y) = \hat{c}(y) + \frac{1}{2}\varepsilon^2(\alpha + 1)$. This implies
$$
\mathcal{L}W(x,y) \leq |x|^{-\alpha}\left ( \mathcal{L}^Y U(y) + \alpha \bar{c}(y) U(y) \right ) .
$$
The rest of the proof follows that of \Cref{thm:ly} by letting $c = \bar{c}$ and $V_\alpha = |x|^{-\alpha}$.
\end{proof}

\begin{assumption}\label{assumption2} Equation \eqref{eq:sp} satisfies the following conditions:
\begin{enumerate}
    \item[(A)] For sufficiently large $|x|$, $x\cdot f(x,y) \leq \hat{c}(y) |x|^2$ for some function $\hat{c}(y)$.
    \item[(B)] There exists a pair of bounded regular domains $\Omega_\text{I}\subset\Omega_\text{O}\subseteq\mathbb{R}^{d_2}$ and an exterior Lyapunov function $\tilde{U}\geq 1$ such that 
    \begin{equation}
          \mathcal{L}^{Y_0}\tilde{U}(y)\leq-p(y)\tilde{U}(y), \quad \forall y\notin\Omega_\text{I},
    \end{equation}
    and
    \begin{equation}
        \mathcal{L}^{Y_0}\tilde{U}(y) + \alpha\bar{c}(y)\tilde{U}(y)\leq q(y), \quad \forall y\in\Omega_\text{I}
    \end{equation}
    for all sufficiently small $\alpha>0$, some positive $C^2$ function $q(y)$, and some uniformly positive function $p(y)$ that satisfies
     \begin{equation}
            p(y) > \zeta \bar{c}(y) \quad \mbox{and} \quad p(y) > \zeta \frac{|\nabla \tilde{U}(y)|}{\tilde{U}(y)}
    \end{equation}
    for some constant $\zeta>0$ uniformly for $y\in\Omega_\text{O}$ and $\bar{c}(y) = \hat{c}(y) + \frac{1}{2}\varepsilon^2\phi(\alpha)$ where $\phi(\alpha)=\alpha-1$ if $\alpha\geq 2$ and $\phi(\alpha)=1+(\alpha-2)/d_1$ if $0<\alpha<2$.
    \item[(C)]  The Feynman--Kac equation
    \begin{equation}
    \label{FKsmall2}
    \begin{cases}
        \begin{aligned}
            &\mathcal{L}^{Y_0}Q + \alpha \bar{c} Q = -q(y)-\gamma, && \text{in }\Omega_\text{I}\\
            &Q = \varphi, && \text{on } \partial\Omega_\text{I}
        \end{aligned}
    \end{cases}
        \end{equation}
    has a positive solution in $C^2(\Omega_\text{I})$ for all $\alpha>0$, where $\varphi\in C^3$ is a positive function defined on $\partial\Omega_\text{I}$ and $\gamma>0$ is a constant.
    \item[(D)] For some $\zeta/3>\alpha>0$, the solution $Q(y)$ to equation \eqref{FKsmall2} has a $C^2$ extension to the whole domain such that 
    \begin{equation}
    \begin{cases}
        \begin{aligned}
            &\mathcal{L}^{Y_0}Q + \alpha\bar{c}Q \leq \frac{1}{2}p(y)\tilde{U}(y), && y\in\Omega_\text{O}\setminus\Omega_\text{I}\\
            &Q = C, && y\notin\Omega_\text{O}
        \end{aligned}
    \end{cases}
        \end{equation}
    for some constant $C$ such that $\tilde{U}\geq C$ on $\Omega_\text{O}^c$.
    \item[(E)] $\lim_{|x|\to \infty}h(x,y)= 0$, where $h(x,y)$ is given in \eqref{approximation}.
\end{enumerate}
\end{assumption}

\begin{cor}
\label{cor2}
    \Cref{assumption2} implies that $|x|^{\alpha}U$ is a Lyapunov function for equation \eqref{eq:sp} provided $|x|$ is sufficiently large, where $U = \tilde{U} + Q$.
\end{cor}
\begin{proof}
    Define $W(x,y)=|x|^{\alpha}U(y)$ and $U=\tilde{U}+Q$, which leads to
    \begin{equation}
        \mathcal{L} W(x,y) = U(y)\mathcal{L}^X  |x|^{\alpha} + |x|^{\alpha}\mathcal{L}^Y U(y).
    \end{equation}
    Direct computations yield
$$
    \mathcal{L}^{X} |x|^{\alpha} = \alpha |x|^{\alpha - 2}\left[x \cdot f(x,y)\right] + \frac{1}{2}\varepsilon^2 \alpha |x|^{\alpha} \left[1+(\alpha-2)\frac{\sum_{i=1}^{d_1}x_i^4}{|x|^4}\right]
$$
Elementary calculation gives
$$
\frac{1}{d_1} \leq \frac{\sum_{i=1}^{d_1}x_i^4}{|x|^4} \leq 1 .
$$
Therefore, we have
$$
\mathcal{L}^{X} |x|^{\alpha}  \leq \alpha |x|^{\alpha}\hat{c}(y) + \frac{1}{2}\varepsilon^2\alpha \phi(\alpha)|x|^{\alpha} := \alpha |x|^\alpha \bar{c}(y)
$$
for
$$
\phi(\alpha) = \left \{\begin{array}{ll}
     \alpha - 1& \quad \alpha \geq 2 \\
     1 + \frac{\alpha - 2}{d_1}& \quad 0 < \alpha  < 2 ,
\end{array}
\right. 
$$
where $\bar{c}(y)=\hat{c}(y)+\frac{1}{2}\varepsilon^2\phi(\alpha)$. The rest of the proof follows that of \Cref{thm:ly} by letting $c=\bar{c}$ and $V_\alpha = |x|^\alpha$.
\end{proof}

\section{Machine-learning-based solver for the Feynman--Kac equation}
\label{sec:numerics}

In this section, we develop a numerical approach to compute the prefactor $Q(y)$, which appears in the Lyapunov function and satisfies the Feynman--Kac equation \eqref{eq:feynman} with Dirichlet boundary conditions. The knowledge of $Q(y)$ is crucial because \Cref{NLS27} only gives the existence of a solution to the Feynman--Kac equation \eqref{eq:feynman}, not its $C^2$ extension to the whole domain. Therefore, for the general system \eqref{SDE2} that is not one-way coupled, we need the numerical solution of the Feynman--Kac equation to show that it admits a $C^2$ extension such that $\mathcal{L}^{Y_0}Q$ is dominated by $\mathcal{L}^{Y_0}\tilde{U}$. 

Our method combines Monte Carlo simulations and the physics-informed neural network (PINN) framework. The Monte Carlo simulations provide approximate values of $Q(y)$ at pre-selected points based on its stochastic representation, while the neural network enforces the governing equation and boundary conditions to obtain a continuous approximation of $Q(y)$ over the domain $\Omega_\text{I}$, as well as its extension to the whole domain.
This method follows a similar strategy for solving high-dimensional Feynman--Kac equations as that developed for high-dimensional Fokker–Planck equations in \cite{zhai2022deep}.

\subsection{Monte Carlo sampling and training data}
We first approximate the solution $Q(y)$ of the Feynman--Kac equation \eqref{eq:feynman} using Monte Carlo simulations based on its stochastic representation \eqref{eq:sto_rep}. Specifically, we evaluate the expectations at a small set of collocation points inside the inner domain $\Omega_\text{I}$, denoted by $\mathcal{Y}^m = \{ y^m_1, \ldots, y^m_{K_1} \} \subset \Omega_\text{I}$. For each $y^m_i$, we simulate multiple trajectories of the SDE and compute the empirical mean, yielding an approximate solution $Q^m_i \approx Q(y^m_i)$, forming $\mathcal{Q}^m = \{ Q^m_1, \ldots, Q^m_{K_1} \}$. 
To impose the Dirichlet boundary condition, we preselect a set of boundary points $\mathcal{Y}^b = \{ y^b_1, \ldots, y^b_{K_2} \} \subset \partial \Omega_\text{I}$. We select boundary points to ensure good coverage of the boundary geometry and pair them with the corresponding boundary conditions $\mathcal{Q}^b = \{\varphi(y^b_1), \ldots, \varphi(y^b_{K_2})\}$. 
In addition, we select a set of interior points inside the inner domain as $\mathcal{Y}^d = \{ y^d_1, \ldots, y^d_{K_3} \} \subset \Omega_\text{I}$, which do not have corresponding $Q$ values. We enforce the PDE constraint by minimizing the residual of the Feynman--Kac equation on these interior points. To enforce the extension condition \eqref{eq:extension}, we generate a set of exterior points $\mathcal{Y}^e=\{ y^e_{1}, \ldots, y^e_{K_4} \}\subset \Omega_\text{O}^c$ and require the network output to approach a prescribed constant $C$ on $\Omega_\text{O}^c$. 
In addition, we consider a set of intermediate points $\mathcal{Y}^f=\{ y^f_{1}, \ldots, y^f_{K_5} \}\subset \Omega_\text{O}\setminus \Omega_\text{I}$, which do not have associated Monte Carlo values, to control the behavior of the extended solution in the region between the inner and outer domains. 

Together, the five datasets $(\mathcal{Y}^m,\mathcal{Q}^m)$, $(\mathcal{Y}^b,\mathcal{Q}^b)$, $\mathcal{Y}^d$, $(\mathcal{Y}^e,C)$, and $\mathcal{Y}^f$ are the training data to solve the Feynman--Kac equation and its extension to the whole domain. Each dataset serves a distinct role: $(\mathcal{Y}^m, \mathcal{Q}^m)$ provides reference values for the solution $Q(y)$ based on Monte Carlo simulations, $(\mathcal{Y}^b,\mathcal{Q}^b)$ enforces the Dirichlet boundary condition, $\mathcal{Y}^d$ ensures that the learned function satisfies the PDE within the inner domain, $(\mathcal{Y}^e,C)$ enforces the constant extension outside $\Omega_\text{O}$, and $\mathcal{Y}^f$ regularizes the extension in $\Omega_\text{O}\setminus\Omega_\text{I}$ by controlling the generator. 
Generating these datasets presents several challenges. Monte Carlo simulations are computationally expensive and subject to sampling noise. Moreover, interior collocation points must provide sufficient domain coverage to accurately constrain the PDE operator. These challenges motivate the design of a training strategy that balances data fidelity, boundary adherence, and differential consistency, which we describe in the next section. 

\subsection{The loss function}
We approximate the solution of the Feynman--Kac equation \eqref{eq:feynman}  using an artificial neural network $Q_{\bm{\theta}}(y)$, parameterized by the parameters ${\bm{\theta}}$. The neural network maps $y$ to the predicted solution $Q_{\bm{\theta}}(y)$, which is trained to satisfy the stochastic, boundary, and differential constraints introduced in the previous subsection. To determine the optimal parameters ${\bm{\theta}}$, we minimize five loss functions corresponding to the Monte Carlo data, boundary conditions, PDE residuals, and extension constraints:
\begin{equation}
\begin{aligned}
L_1({\bm\theta}) &= \frac{1}{K_1} \sum_{i = 1}^{K_1} \left[ Q_{\bm \theta}(y^m_i) - Q^m_i \right]^2,  \\
L_2({\bm\theta}) &= \frac{1}{K_2} \sum_{i = 1}^{K_2} \left[ Q_{\bm \theta}(y^b_i) - Q^b_i \right]^2, \\
L_3^{Y_0}({\bm\theta}) &= \frac{1}{K_3} \sum_{i = 1}^{K_3} \left[ \mathcal{L}Q_{\bm \theta}(y^d_i) + \alpha c(y^d_i)Q_{\bm \theta}(y^d_i) + q(y^d_i) +\gamma \right]^2, \\
L_4({\bm\theta}) &= \frac{1}{K_4} \sum_{i = 1}^{K_4} \left[ Q_{\bm \theta}(y^e_i) - C \right]^2,  \\
L_5({\bm\theta}) &= \frac{1}{K_5} \sum_{i = 1}^{K_5} \left[ \max\left(0,\mathcal{L}^{Y_0}Q_{\bm \theta}(y^f_i) + \alpha c(y_i^f)Q_{\bm{\theta}}(y_i^f) - \frac{1}{2} p(y^f_i)\tilde{U}(y^f_i) \right)\right]^2.
\end{aligned}
\end{equation}
The first loss function $L_1$ enforces consistency between the neural network predictions and the Monte Carlo data, while $L_2$ ensures that the Dirichlet boundary conditions are satisfied. The third loss function $L_3$ measures the residual of the Feynman--Kac equation inside the domain, enforcing the underlying PDE constraint. The fourth loss function $L_4$ enforces the constant extension $Q_{\bm\theta}= C$ on $\Omega_\text{O}^c$. The fifth loss function $L_5$ controls the growth of the extended solution in $\Omega_\text{O}\setminus\Omega_\text{I}$ by constraining the generator through the exterior Lyapunov function $\tilde{U}$.

It is challenging to balance 5 loss functions during neural network training. In particular, we must prevent the solver from collapsing to the trivial solution $Q = 0$. We take the following steps to train the neural network. {\bf (1) Pre-training}: we use three $L^2$ losses $L_1, L_2$, and $L_4$ to train the neural network for a few epochs. Since the operator $\mathcal{L}$ is not involved, this step is very fast. Although not accurate, optimizing $L^2$ losses will bring the neural network to a suitable baseline. {\bf (2) Alternating Adam}: After that, we adopt the training strategy discussed in \cite{zhai2022deep} by minimizing loss functions $L_1, L_2, L_3,$ and $L_4$ alternately using the Adam optimizer to balance their contributions, as the magnitude of each term varies across training phases. The Adam optimizer's approximate scale invariance naturally maintains stability and prevents one loss component from dominating the optimization. {\bf (3) Fine-tuning}: After the second step, all loss functions become small enough. Then we add the regularization term $L_5$ into the training, and minimize $L := L_1 + L_2 + L_3 + L_4 + L_5$ directly.

In our numerical experiments, we employ a small feedforward neural network for simplicity, since the primary focus of this paper is the Feynman--Kac--Lyapunov framework. The collocation points $\mathcal{Y}^d$ are uniformly sampled within the domain. More sophisticated neural architectures and adaptive sampling strategies can be found in \cite{cai2021physics, meng2020composite, jagtap2020adaptive, wang2021understanding, lu2021physics, gao2021phygeonet}, and guidelines for selecting suitable collocation points are discussed in \cite{wang2024deep}. 
We present the numerical results of the proposed method for a 1D toy example in \Cref{sec:1dtoy} and for a three-mode cascade example in \Cref{sec: 4dcascade}.

\subsection{Improvement of estimates}\Cref{general_assumption}(E) requires the driving force $g(x,y)$ to approach $g_0(y)$ as $x$ goes to its asymptotic limit. This requirement is mainly due to the lack of knowledge of the prefactor $U$ itself. If the prefactor $U$ is known after numerically solving $Q$, we can relax the conditions in \Cref{assumption1}(E) and \Cref{assumption2}(E) and replace them with the following conditions.

\quad

\noindent \Cref{assumption1}($\tilde{\text{E}}$): For sufficiently small $|x|$, we have
\begin{equation}
h(x,y) |\nabla U(y)| \leq \frac{\gamma}{2}, \quad \forall y\in\Omega_\text{I} ,
\end{equation}
\begin{equation}
h(x,y) |\nabla U(y)| \leq \frac{1}{24}p(y) \tilde{U}, \quad \forall y\in\Omega_\text{O}\setminus \Omega_\text{I} ,
\end{equation}
and 
\begin{equation}
    h(x,y) \leq \frac{1}{12}\zeta , \quad \forall y\in\Omega_\text{O}^c .
\end{equation}
\Cref{assumption2}($\tilde{\text{E}}$) can be given analogously. The relaxed conditions make the Feynman--Kac--Lyapunov method more applicable to a more general class of problems, but require the estimation of the upper bound of $\|\nabla U\|_\infty$, which is not an easy task because $U = \tilde{U} + Q$, where $Q$ solves the Feynman--Kac equation. It is desirable to numerically compute the prefactor function $Q(y)$ to use the relaxed conditions \Cref{assumption1}($\tilde{\text{E}}$) and \Cref{assumption2}($\tilde{\text{E}}$).

\section{Examples}\label{sec:example}

We provide three examples of the driving-driven SDEs to illustrate our Feynman--Kac--Lyapunov method to construct the Lyapunov functions. We use both analytical and numerical approaches to show the existence of a positive solution for the Feynman--Kac equation in these cases. 

\subsection{Mobility-based SIR Model} 
Consider an epidemiological compartmental model, where individuals in a population transition between the three states (susceptible, infected, and recovered) based on an infection rate $\beta$, a recovery rate $\gamma$, and their mobility factor $m\in(0,1)$. The mobility-based SIR model was originally proposed and studied in \cite{mSIR2024}. 
We divide the population into $N$ groups based on their mobility $m_n, n=1,\ldots,N$. Let $S_n,I_n,R_n$ be the populations of susceptible, infected, and recovered compartments with the $n$-th mobility group, respectively. We assume the mass-action law, so the infection rate between $S_n$ and $I_{\tilde{n}}$ is proportional to $m_n m_{\tilde{n}}$. The populations of three compartments satisfy
\begin{align*}\begin{cases}
    \mathrm{d}S_n &= -\beta\sum_{\tilde{n}=1}^N m_n m_{\tilde{n}} S_nI_{\tilde{n}}\,\mathrm{d}t\\
    \mathrm{d}I_n &= \beta\sum_{\tilde{n}=1}^N m_n m_{\tilde{n}}S_nI_{\tilde{n}}\,\mathrm{d}t - \gamma I_n\,\mathrm{d}t\\
    \mathrm{d}R_n &= \gamma I_n\,\mathrm{d}t.
    \end{cases}
\end{align*}
As the time-dependent mobility population is not accessible in practice, we replace the population mobility with the empirical distribution of $M$ independent agents. Let $Y_t = (Y_{t}^{(1)}, Y_{t}^{(2)}, \ldots,Y_{t}^{(M)})$ be the mobility of $M$ agents, satisfying
\begin{equation}
\label{eq:SIR-Y}
    \mathrm{d}Y_{t}=g(Y_{t})\,\mathrm{d}t + \sigma\,\mathrm{d}B_t,
\end{equation}
where $B_t$ is the Wiener process in $\mathbb{R}^M$. 
When the agents are sampled from a large population, the mobility of group $n$ is determined by the empirical mobility distribution of observed agents, which replaces $m_n$ with $m(n, Y_t)$, where $m(n, \cdot)$ is a bounded and nonnegative measurable function for each $n = 1,\ldots, N$. We then obtain the agent-driven mobility-SIR model
\begin{align}\label{eq:mSIR}\begin{cases}
    \mathrm{d}S_n = -\beta\sum_{\tilde{n}=1}^N m(n,Y_t)m(\tilde{n}, Y_t)S_nI_{\tilde{n}}\,\mathrm{d}t\\
    \mathrm{d}I_n = \beta\sum_{\tilde{n}=1}^N m(n,Y_t)m(\tilde{n}, Y_t)S_nI_{\tilde{n}}\,\mathrm{d}t - \gamma I_n\,\mathrm{d}t\\
    \mathrm{d}R_n = \gamma I_n\,\mathrm{d}t\\
    \mathrm{d}Y_t=g(Y_t)\,\mathrm{d}t + \sigma\,\mathrm{d}B_t,
\end{cases}\end{align}
where $m(n,Y_t)$ is a function that gives the mobility of group $n$ based on the observation of agents. There are many choices for defining the mobility function $m(n, Y_t)$. For example, let $m(n, Y_t)$ be the mean mobility of the $n$-th percentile of agents. More precisely, we rank the mobility of agents such that 
\begin{equation}
Y_{t}^{(k_1)} \ge Y_{t}^{(k_2)} \ge \ldots \ge Y_{t}^{(k_M)}.
\end{equation}
We let $M=n\triangle k$ and define
\begin{equation}
m(n,Y_t) = \frac{1}{\triangle k}\sum_{i = (n-1)\triangle k+1}^{n\triangle k} Y_{t}^{(k_i)} .
\end{equation}

Equation \eqref{eq:mSIR} is a driving-driven system. Assume the following conditions hold. 

\begin{assumption}\label{assumption3}
    Let $\mathbf{I}$ be the column vector of the infected population. Assume equation \eqref{eq:mSIR} satisfies the following conditions.
\begin{enumerate}
    \item[(A)] The mobility dynamics equation \eqref{eq:SIR-Y} is geometrically ergodic with a unique invariant probability measure $\pi$ such that $\pi(V^2_y)<\infty$ for a Lyapunov function $V_y$.
    \item[(B)] The invariant probability measure $\pi$ satisfies
\begin{equation}
    \beta\,\mathbb{E}_{\pi}\left[J\right]-\gamma <0 \quad \text{and} \quad \mathrm{Var}_{\pi}\left(J\right)>0 ,
\end{equation}
where 
\begin{equation}
    \label{eq:J}
    J = \frac{1}{2}\left ( \sum_{n = 1}^N m^2(n, Y_t)S_n(0) + \sqrt{\sum_{n = 1}^N m^2(n, Y_t)S^2_n(0) \sum_{n = 1}^N m^2(n, Y_t)}\right ) .
\end{equation}
\end{enumerate}
\end{assumption} 
\begin{cor}
    \label{ex1}
    Assume $R_n(0) = 0$  for all $n = 1, \ldots, N$. If Assumption \ref{assumption3} holds, there exists a strictly positive function $U(Y_t)$ and a sufficiently small constant $\alpha > 0$ such that $|\mathbf{I}|^\alpha U$ is a Lyapunov function for the SDE \eqref{eq:mSIR}.
\end{cor} 
\begin{proof}
Let
\begin{equation}
    \gamma - \beta\,\mathbb{E}_{\pi}\left[J\right] := \eta > 0 .
\end{equation}
Define $W(\mathbf{I},Y_t)=|\mathbf{I}|^{\alpha}U(Y_t)$, where the function $U$ will be decided later in the proof. Computing $\mathcal{L}W$, we have
\begin{equation}
\begin{aligned}
    \mathcal{L}W & = U(Y_t)\mathcal{L}^{\mathbf{I}}|\mathbf{I}|^{\alpha} + |\mathbf{I}|^{\alpha}\mathcal{L}^YU(Y_t)\\
    &=\alpha|\mathbf{I}|^{\alpha}\frac{\mathbf{I}^TA(\mathbf{S},Y_t)\mathbf{I}}{\mathbf{I}^T\mathbf{I}}U(Y_t)+\mathcal{L}^YU(Y_t)|\mathbf{I}|^{\alpha}
\end{aligned}
\end{equation}
where 
\begin{equation}
    A(\mathbf{S},Y_t)=\beta
    \begin{pmatrix}
        m_1S_1\\
        m_2S_2\\
        \vdots\\
        m_NS_N
    \end{pmatrix}
    \begin{pmatrix}
        m_1 & m_2 & \ldots & m_N
    \end{pmatrix} -
    \gamma \mathrm{Id} 
\end{equation}
where $m_n=m(n,Y_t)$.

Notice that $A(\mathbf{S}_0,Y_t)$ is the difference of a rank-one matrix and $\gamma \mathrm{Id}$. The Rayleigh quotient of the rank-one matrix equals the Rayleigh quotient of its symmetric part, which is $(\mathbf{u}\mathbf{m}^T + \mathbf{m}\mathbf{u}^T)/2$ for $\mathbf{m} = (m_1, \ldots, m_n)$ and $\mathbf{u} = (m_1 S_1(0), \ldots, m_n S_n(0))$. Further noticing that all entries of $I$ and $A$ are nonnegative, we have
\[
\frac{\mathbf{I}^TA(\mathbf{S},Y_t)\mathbf{I}}{\mathbf{I}^T\mathbf{I}}  \leq \frac{\mathbf{I}^TA_0(\mathbf{S},Y_t)\mathbf{I}}{\mathbf{I}^T\mathbf{I}} \leq \frac{1}{2}(|\mathbf{u}||\mathbf{m}| + \mathbf{u}^T\mathbf{m}) 
\]
because $S_n(t) \leq S_n(0)$. This gives
\begin{equation}
    \frac{\mathbf{I}^TA(\mathbf{S},Y_t)\mathbf{I}}{\mathbf{I}^T\mathbf{I}} \leq \beta J - \gamma  .
\end{equation}
Now let $V_{\alpha} = | \mathbf{I}|^{\alpha}$, we have
\begin{equation}
    \mathcal{L}^I V_{\alpha} \leq \alpha \left ( \beta J -\gamma  \right) V_{\alpha}
\end{equation}
Note that $ c(Y_t):=\beta J-\gamma $
is uniformly bounded and has nontrivial asymptotic variance with respect to $\pi$ by \Cref{assumption3}.
It remains to show that the $Y_t$-dynamics is geometrically ergodic and has a Lyapunov function $V$ satisfying $\pi(V^2) < \infty$, which is given in \Cref{assumption3}. 

Therefore, by \Cref{thm:oneside}, when $\alpha$ is sufficiently small, the Feynman--Kac operator $\mathcal{L}^Y + \alpha c$ admits a strictly positive eigenfunction $U_\alpha$ such that $|I|^\alpha U_\alpha$ is the Lyapunov function of equation \eqref{eq:mSIR}. This completes the proof. 
\end{proof}

Let $\mathcal{L}$ be the generator of equation \eqref{eq:mSIR}. Under \Cref{assumption3}, because $\mathcal{L}|\mathbf{I}|^\alpha U < 0$ for all $|\mathbf{I}|\neq 0$, the total number of infected individuals decreases to $0$ as $t\rightarrow \infty$. In other words, the disease-free equilibrium is asymptotically stable. 

\subsection{The 1D oscillator}\label{sec:1dtoy}
We use the following toy example to illustrate our numerical method. Consider the following system where $X_t$ is driven by a 1D oscillator $\theta_t$ on $\mathbb{S}_1$.
\begin{equation}
\begin{aligned}
    \label{1D}
    \mathrm{d}X_t &= (\sin \theta_t + \lambda) X_t \,\mathrm{d}t\\
    \mathrm{d}\theta_t & = -(2 \cos \theta_t+1) \,\mathrm{d}t + \varepsilon \,\mathrm{d}B_t,
\end{aligned}
\end{equation}
where $\lambda$ is a positive constant that we discuss later. 

Notice that the $\theta$-dynamics concentrates at the stable equilibrium $\theta^* = 4\pi/3$. For $\lambda < \sqrt{3}/2$, we have $\sin \theta^* + \lambda < 0$ and thus $X_t$ should be stable. Suppose the Lyapunov function $W(x, \theta) = |x|^\alpha U(\theta)$ has a product form. Since the domain of $\theta$ is compact, there is no need for an extension or an exterior Lyapunov function. Applying the infinitesimal generator of equation \eqref{1D} to $W$, we have
$$
\mathcal{L} W = |x|^\alpha \left [ -(2 \cos \theta+1)U_\theta + \frac{1}{2}\varepsilon^2 U_{\theta\theta} + \alpha (\sin \theta + \lambda)U\right] .
$$
Provided that $U$ satisfies the 1D Feynman--Kac equation
\begin{equation}
    \label{1D-FK}
    -(2 \cos \theta+1)U_\theta + \frac{1}{2}\varepsilon^2 U_{\theta\theta} + \alpha (\sin \theta + \lambda)U = -1,
\end{equation}
then $W$ is a Lyapunov function satisfying 
$$
\mathcal{L} W = - |x|^\alpha.
$$
Theoretically, the first half of the proof in \Cref{NLS27} still holds, and we know that the operator 
\begin{equation}
\mathcal{L}^{\theta} := -(2 \cos \theta + 1) \partial_\theta + \frac{1}{2}\varepsilon^2 \partial^2_{\theta} + \alpha (\sin \theta + \lambda)
\end{equation}
admits a strictly negative principal eigenvalue for all sufficiently small $\alpha$, if $\mathbb{E}_\pi[\sin \theta + \lambda] < 0$, where $\pi$ is the invariant probability measure of the $\theta$-variable. Then by the Fredholm alternative, $\mathcal{L}^{\theta}$ is invertible, so $\mathcal{L}^{\theta} U = -1$ admits a $C^2$ solution. And by equation \eqref{eq:sto_rep}, $U$ is strictly positive. However, the caveat is that (1) the eigenvalue argument only works for {\it sufficiently small} $\alpha > 0$ and (2) we can only prove $\mathbb{E}_\pi[\sin \theta + \lambda] < 0$ for {\it sufficiently small $\varepsilon$}.

\medskip

Here we compute the prefactor function $U$ using the machine-learning approach and show how model parameters affect the properties of $U$. Our numerical results suggest that the positive solution of the Feynman--Kac equation exists for a much wider range of parameters. 

In this numerical example, we choose $\alpha = 0.5$ and $\lambda = 0.1$. We use a small feedforward artificial neural network to represent the solutions of equation \eqref{1D-FK}. The Monte Carlo training set $(\mathcal{Y}^m, \mathcal{U}^m)$ contains $50$ pairs of data points, while each $u_i^m$ is obtained by averaging over $1000$ samples of the stochastic representation in equation \eqref{eq:sto_rep}. For the PDE constraint dataset, we randomly sample $1000$ points uniformly from $[-\pi, \pi]$ to form $\mathcal{Y}^d$ and minimize the PDE residual over these points.

\begin{figure}[htp]
    \centering
    \includegraphics[width=0.8\linewidth]{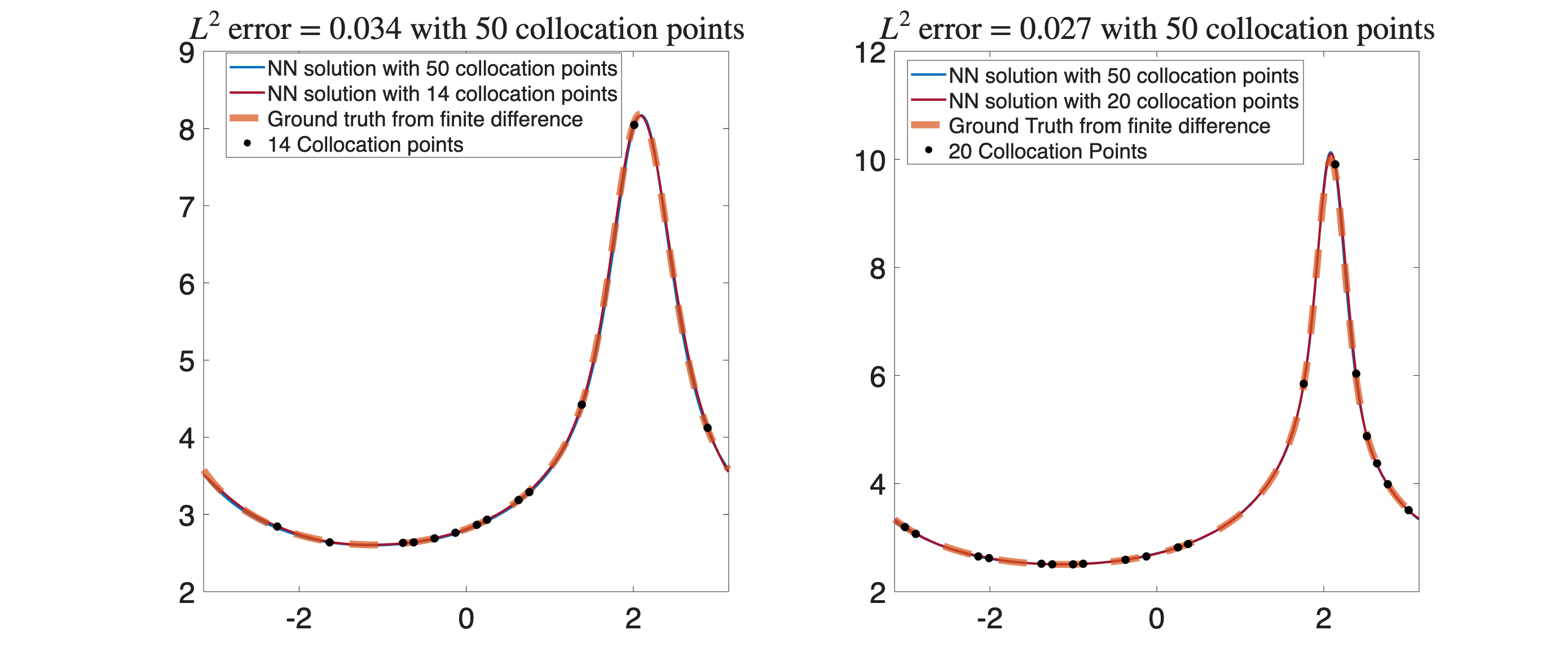}
    \caption{Numerical solutions to equation \eqref{1D-FK} for (left) $\varepsilon = 0.5$ and (right) $\varepsilon = 0.25$, using a neural network approach and the finite difference method. }
    \label{fig1}
\end{figure}

In \Cref{fig1}, we show the numerical solution obtained from the artificial neural network and the ground truth computed by a finite difference scheme. The left and right panels are for different noise levels, $\varepsilon = 0.5$ and $\varepsilon = 0.25$, respectively. 
When the noise is weak, the stochastic dynamics stay at the unstable equilibrium for a longer time before entering the stable equilibrium, giving a larger value of $u$ in the small neighborhood of the unstable equilibrium $2\pi/3$. Therefore, the solution is much more concentrated near the unstable equilibrium when $\varepsilon = 0.25$. 

\begin{figure}[htp]
    \centering
    \includegraphics[width=\linewidth]{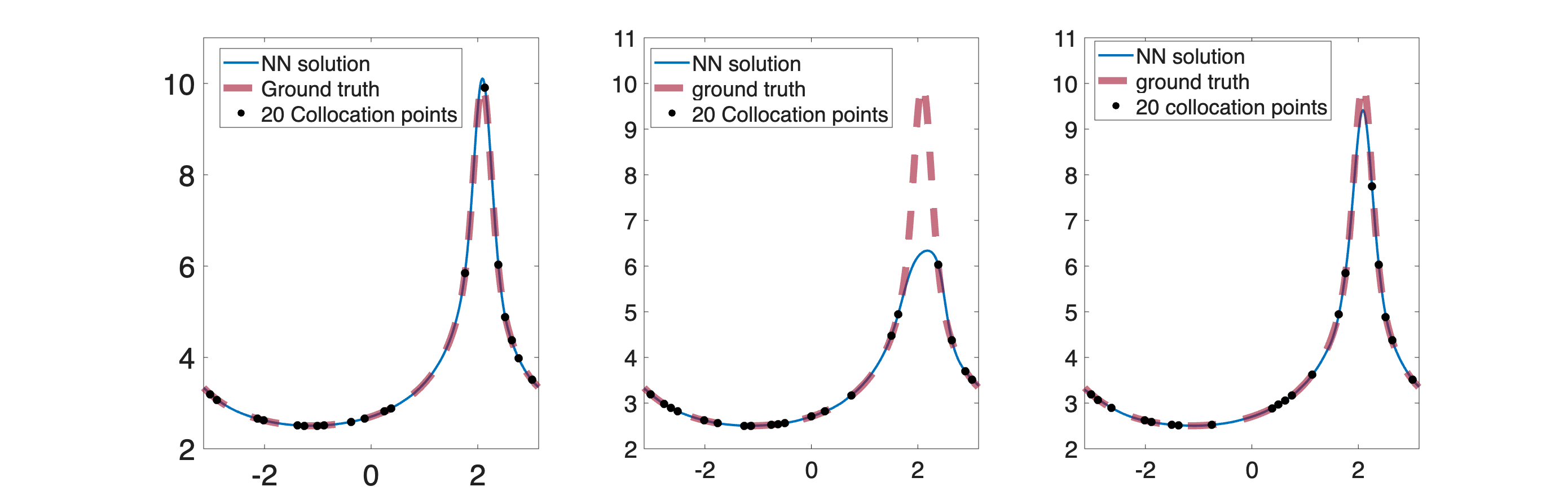}
    \caption{Numerical solutions to equation \eqref{1D-FK} using three different sets of collocation points. }
    \label{fig2}
\end{figure}

Then we randomly choose $14$ collocation points when $\varepsilon = 0.5$ and $20$ collocation points when $\varepsilon = 0.25$ to demonstrate the performance of the neural network solver with fewer collocation points. Note that we take more collocation points when $\varepsilon = 0.25$ because the solution is more concentrated. As seen in \Cref{fig1}, the solution with $\varepsilon = 0.5$ is very accurate, while there are some visible errors in the solution with $\varepsilon = 0.25$ near the unstable equilibrium. Generally speaking, an artificial neural network has higher approximation error when approximating functions with large derivatives. To further demonstrate the effect of limited sample size, we select another two subsets of collocation points and run the neural network solver again. One subset of collocation points does not contain the peak near $4\pi/3$ at all, while the other only has one collocation point near the top of the peak. As seen in \Cref{fig2}, the quality of those two solutions is inferior to the first one, especially when the peak is not well-represented. This highlights the need to select a suitable collocation point set with suitable concentration at regions where the solution has dramatic changes. 

\subsection{The three-mode energy cascade model} \label{sec: 4dcascade}
In this subsection, we study a three-mode energy cascade model derived from the nonlinear Schr\"odinger equation. The model in \cite{colliander2010transfer} describes an energy cascade model with nearest-neighbor interactions and is derived from the cubic nonlinear Schr\"odinger equation to demonstrate the energy cascade dynamics that transfers energy from low modes to high modes. We focus on the three-mode version of this system and rewrite it in action–angle coordinates. The variables $I_1$, $I_2$, and $I_3$ are the action variables, and $\theta_1$ and $\theta_3$ are the phase differences between neighboring modes. We refer to \cite{hani2025nonequilibrium} for detailed discussions of the energy cascade model, its role in mathematical physics, and the transformation to action–angle coordinates.

We consider external noise and damping in this three-mode energy cascade model and obtain the following nonequilibrium energy cascade model, where
\begin{equation}
  \label{cascade}
  \begin{aligned}
  \mathrm{d}I_{1}& = [2 I_{1} I_{2} (\sin \theta_{1} -
  \rho) + \rho(T_{1} - I_{1}^2) ] \,\mathrm{d}t + \sqrt{\rho T_{1} I_{1}/2} \,\mathrm{d}B_1 \\
  \mathrm{d}I_{2} &= - 2 I_{2} (I_{1} \sin \theta_{1} + I_{3} \sin \theta_{3}) \,\mathrm{d}t\\
  \mathrm{d}I_{3} &= [2 I_{2}I_{3}(\sin \theta_{3} - \rho) + \rho(T_{3} - I_{3}^2)] \,\mathrm{d}t + \sqrt{\rho T_{3}I_{3}/2} \,\mathrm{d}B_2\\
  \mathrm{d} \theta_{1}& =  \left[I_{2}(1 + 2 \cos \theta_{1}) - I_{1} (1 + 2 \cos \theta_{1}) - 2 I_{3} \cos \theta_{3}\right] \,\mathrm{d}t + \sqrt{\rho}  \,\mathrm{d}B_3\\
  \mathrm{d} \theta_{3} &= \left[I_{2} (1 + 2 \cos \theta_{3}) - I_{3} (1 + 2 \cos \theta_{3}) - 2 I_{1}\cos \theta_{1}\right] \,\mathrm{d}t + \sqrt{\rho}\,\mathrm{d}B_4,
\end{aligned}
\end{equation}
where $0<\rho<1$ is a model parameter. It is generally difficult to establish the stochastic stability of the nonequilibrium energy cascade model. In \cite{hani2025nonequilibrium}, the authors prove stochastic stability for a modified version of equation \eqref{cascade} by introducing additional mechanisms, such as enhanced damping terms, to make a rigorous proof possible. The goal of this subsection is to illustrate how to compute a prefactor function of the Lyapunov function in high-dimensional systems and how the numerically obtained prefactor can improve stability estimates for such complex systems.

One major challenge in the proof of stochastic stability of equation \eqref{cascade} is to show that once $I_2$ reaches $0$, it is not trapped near $0$, which leads to the breakdown of the energy cascade model.
Motivated by this challenge, we view equation \eqref{cascade} as a perturbation of the driving-driven dynamics around $I_2=0$.
In particular, we adopt the notation in \eqref{SDE2} and \Cref{assumption1} and define $X = I_2$, $Y = (I_1, I_3, \theta_1, \theta_3)$, and $c = 2(I_1 \sin\theta_1 + I_3 \sin \theta_3)$. 
In addition, we define the function $g_0$ as
\begin{equation}
g_0(Y) = \begin{bmatrix}
    \rho (T_1 - I_1^2)\\
    \rho (T_3 - I_3^2)\\
    -I_1 (1 + 2 \cos \theta_1) - 2 I_3 \cos \theta_3\\
    -I_3 (1 + 2 \cos \theta_3) - 2 I_1 \cos \theta_1
\end{bmatrix}.
\end{equation}
Therefore, the difference function between $g$ and $g_0$ has the form
\begin{equation}
g(X,Y) - g_0(Y) = \begin{bmatrix}
    2I_1 I_2 (\sin \theta_1 - \rho)\\
    2I_2 I_3 (\sin \theta_3 - \rho)\\
    I_2(1 + 2 \cos \theta_1)\\
    I_2(1 + 2 \cos \theta_3)
\end{bmatrix}.
\end{equation}
Considering the worst case of angles when $\sin \theta_i - \rho = -1 - \rho$ and $1 + 2 \cos \theta_i = 3$ for $i = 1,3$, we get the bound 
\begin{equation}
h=\|g - g_0\|_2 \leq  I_2 \sqrt{4(1+\rho)^2( I_1^2 + I_3^2) + 18} ,
\end{equation}
which approaches $0$ as $I_2 \rightarrow 0$. Note that to avoid confusion, in this section, we explicitly state which vector norm is used. 

It is not difficult to see that when $I_1 + I_3$ is large, $\tilde{U} = 3( I_1 + I_3)$ and $p = 0.4\rho(I_1 + I_3)$ satisfy \Cref{assumption1}(B). In particular, when $I_1 + I_3$ is sufficiently large, it satisfies that 
\begin{equation}
\mathcal{L}^{Y_0}\tilde{U} = 3\rho[(T_1 + T_3) - I_1^2 - I_3^2 ]\leq 3\rho [(T_1 + T_3) - \frac{1}{2}(I_1 + I_3)^2 ] \leq - 0.4 \rho (I_1 + I_3) \tilde{U} .
\end{equation}
Noting that $c \leq 2(I_1 + I_3)$, we have
\begin{equation}
    p(y) = 0.4\rho(I_1+I_3) \geq 2\zeta(I_1+I_3)\geq \zeta c(y)
\end{equation}
\begin{equation}
    p(y)-\zeta\frac{\|\nabla \tilde{U}(y)\|_2}{\tilde{U}(y)} = 0.4 \rho(I_1+I_3)-\zeta\frac{\sqrt{2}}{I_1+I_3}>0
\end{equation}
provided $\zeta \leq 0.2 \rho$, and $\tilde{U} = 3(I_1 + I_3)$ is sufficiently large. Therefore, $\tilde{U}$ is a suitable exterior Lyapunov function. In addition, considering the worst case of $\sin \theta_1$ and $\sin \theta_3$, we have
\begin{align*}
    \mathcal{L}^{Y_0} \tilde{U} + \alpha c \tilde{U} &\leq 3\rho [ (T_1 + T_3) - I_1^2 - I_3^2 + 2 \alpha (I_1 + I_3)^2] \\
    &\leq 3\rho [ (T_1 + T_3) - I_1^2 - I_3^2 + 4 \alpha (I_1^2 + I_3^2) ] \\
    &= 3\rho [(T_1 + T_3) - (1 - 4 \alpha)(I_1^2 + I_3^2) ],
\end{align*}
meaning that $\mathcal{L}^{Y_0} \tilde{U} + \alpha c \tilde{U}$ is negative for large enough $I_1 + I_3$. Letting $\Omega_\text{I} = [0, 12]^2 \times [-\pi, \pi]^2$ and $\Omega_\text{O} = [0, 16]^2 \times [-\pi, \pi]^2$, it is easy to see that $\tilde{U} = 3(I_1 + I_3)$ satisfies conditions (A) and (B) of \Cref{general_assumption} on $\Omega_\text{I}^c$.

\begin{figure}
    \centering
    \includegraphics[width=\linewidth]{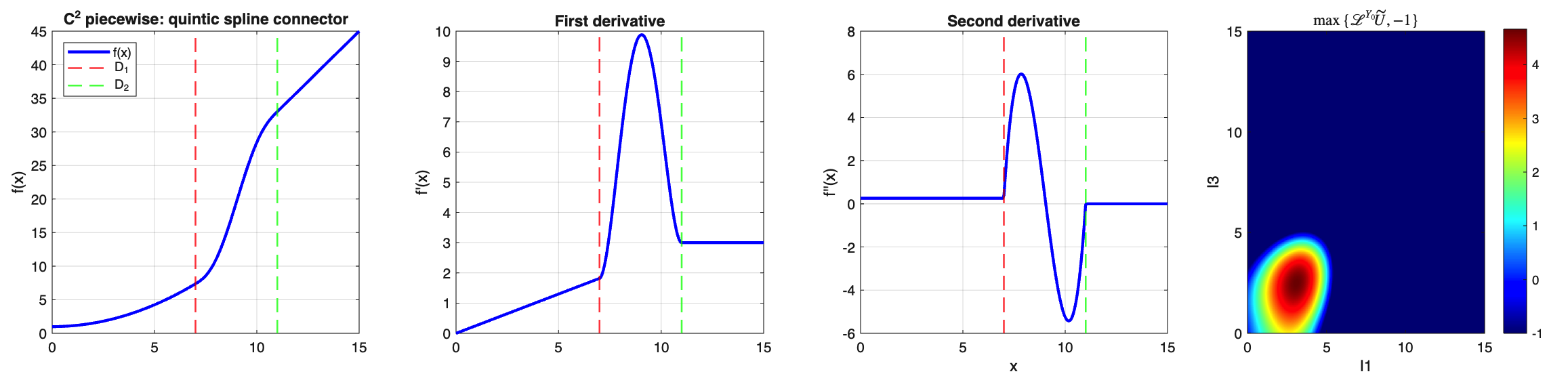}
    \caption{Left to Right: The quintic interpolation of $u(x)$, first derivative of $u$, second derivative of $u$, $\max\{\mathcal{L}^{Y_0} \tilde{U}, -1\}$.}
    \label{fig:tildeU}
\end{figure}

To facilitate computation, we avoid $\mathcal{L}^{Y_0} \tilde{U} + \alpha c \tilde{U}$ being too large within $\Omega_\text{I}$, as the value of the solution to the Feynman--Kac equation scales linearly with $q(y)$. Therefore, we modify $\tilde{U}$ for small $I_1$ and $I_3$ to reduce the value of $\mathcal{L}^{Y_0} \tilde{U}$ by defining
$$
u(x) = \left \{\begin{array}{lll}
     1 + 0.13x^2 & \mbox{if} & 0\leq x < 7 \\
     \mbox{Quintic interpolation} & \mbox{if} & 7 \leq x < 11\\
     3x & \mbox{if } & x \geq 11
\end{array}
\right.
$$
and $\tilde{U} = u(I_1) + u(I_3)$. In order to make $u(x)$ a $C^2$ function, there are $6$ constraining equations at $x = 7$ and $x = 11$. Therefore, there exists a unique quintic interpolation that can be easily solved numerically making $\tilde{U}$ a $C^2$ function. Since $\tilde{U}$ only depends on $(I_1, I_3)$, we can then check the value of $\mathcal{L}^{Y_0} \tilde{U} + \alpha c \tilde{U} $. When $\alpha = 0.025$, $T_1 = 2$, $T_3 = 10$, as demonstrated in \Cref{fig:tildeU}, the maximum value of $\mathcal{L}^{Y_0} \tilde{U} + \alpha c \tilde{U} $ is approximately $4.66$, a lot smaller than $3\rho(T_1 + T_3) = 36\rho$ before modification. In addition, outside of $\Omega_\text{O}$ we have $\tilde{U} \geq 48$. Therefore, we let $C = 48$ be the value of $Q$ outside $\Omega_\text{O}$. 

By \Cref{cor1}, there exists a Lyapunov function $W = I_2^{-\alpha}U$ for all sufficiently small $I_2$, where the prefactor $U=\tilde{U}+Q$, provided $Q$ is the strictly positive solution to equation \eqref{FKsmall} in a bounded domain $\Omega_\text{I}$ and can be extended to the whole domain under \Cref{assumption1}(D). In \cite{hani2025nonequilibrium}, after a careful and tedious calculation of the invariant probability measure of the $Y$-dynamics, we proved the existence of a solution to equation \eqref{FKsmall} when $0 < \alpha \ll 1$, $T_1 \ll T_3$ and $\rho \ll T_3$. In this section, we will numerically compute the Lyapunov prefactor component $Q$, both relaxing theoretically proved results and demonstrating the algorithm in \Cref{sec:numerics} for higher-dimensional examples. 

\begin{figure}[htp]
    \centering
    \includegraphics[width=0.45\linewidth]{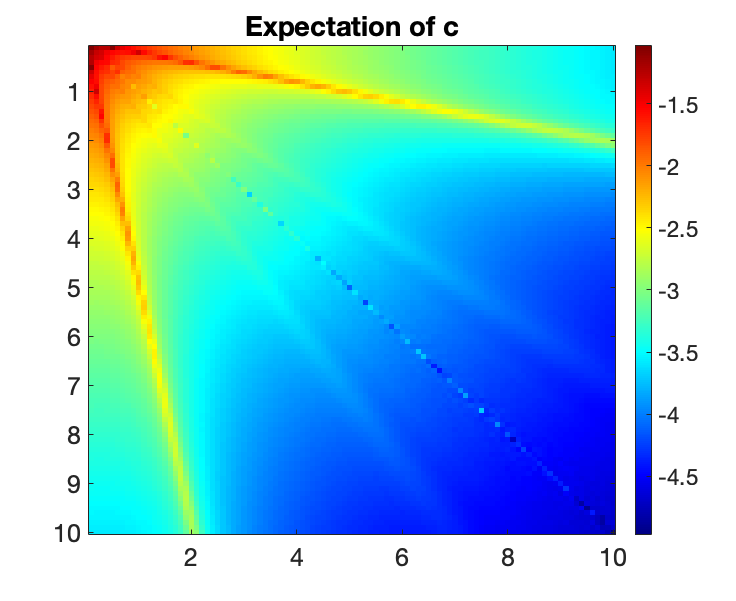}
    \caption{Expectation $\pi^Y(c)$ for $100\times 100$ $(T_1, T_3)$ pairs in $[0,  10]^2$.}
    \label{fig3}
\end{figure}

We first check the sign of the observable $\pi^Y(c)$, where $\pi^Y$ is the invariant probability measure of the $Y$-dynamics. The rigorous proof requires $T_1 \ll T_3$ to make $I_3 \gg I_1$ with dominantly high probability. Then one can show that $\pi^Y(c) < 0$ because $\theta_3$ concentrates at $4\pi/3$ as $I_3 \gg I_1$. However, our numerical results show that $\pi^Y(c) < 0$ regardless of the values of $ T_1$ and $T_3$. We tested combinations of $T_1$ and $T_3$ on a $100\times 100$ grid in $[0,  10]^2$. For each $(T_1, T_3)$ pair, we use Monte Carlo simulation to estimate $\pi^Y(c)$ by integrating $c$ over a long trajectory up to $T = 10^5$. The result is shown in \Cref{fig3}. We can see that $c$ always has negative expectation with respect to $\pi^Y$. This highlights the gap between what actually happens and what can be rigorously proved in numerous applications. 

We then compute the solution $Q$ of the Feynman--Kac equation \eqref{FKsmall} with parameters {$T_1=2$, $T_3=10$, $\rho = 0.4$, $q = 4.67$, $\gamma = 0.13$, and $\alpha = 0.025$}. Recall that the maximal value of $\mathcal{L}^{Y_0} \tilde{U} + \alpha c \tilde{U}$ is approximately $4.66$, satisfying \Cref{assumption1}(B). The computational domain is {$\Omega_\text{I} = [0, 12]^2 \times [-\pi, \pi]^2$}. With this parameter set, we have 
$$
    \mathcal{L}^{Y_0} \tilde{U} = 3\rho\left [ 12 - I_1^2 - 144\right ] \leq - 1.2 \rho (I_1 + 12)^2
$$
for all $I_1 \in [0, 12]$. The case of the other edge of the boundary, which is $(12, I_3)$ for $I_3 \in [0, 12]$, is analogous. Then we impose the natural boundary condition, which coincides with the solution to the Feynman--Kac equation 
$$
    \mathcal{L}^{Y_0} u + \alpha c u = - q - \gamma \quad \mbox{ in } \quad \mathbb{R}^2_+ \times \mathbb{T}^2
$$
at $\partial \Omega$. This allows us to run the stochastic representation in equation \eqref{eq:sto_rep} with $\tau \rightarrow \infty$ in the whole space $\mathbb{R}^2_+ \times \mathbb{T}^2$. The simulation stops when the exponent $\int_0^\tau c_0(X_s) \mathrm{d}s$ is less than $-10$, indicating that the remaining tail is sufficiently small. The boundary segments ${ I_1 = 0 }$ and ${ I_3 = 0 }$ contain non-regular points, meaning that trajectories starting from these boundaries remain inside the domain for some time and do not leave immediately under equation \eqref{cascade}. The Feynman--Kac equation is not uniformly elliptic on $\Omega_\text{I}$. However, we can still verify the existence of a numerical solution that satisfies the requirement of \Cref{thm:ly} and ensures the construction of the Lyapunov function. In addition to solving $Q$ numerically using an artificial neural network, we also use the same artificial neural network to extend it into $\Omega_\text{O} = [0, 16]^2\times [-\pi, \pi]^2$ with $Q|_{\Omega_\text{O}^c} = C = 48$, which is less than the value of $\tilde{U}$ at the boundary of $\Omega_\text{O}$. Hence \Cref{assumption1}(D) is satisfied. 

For the Monte Carlo simulation training data $(\mathcal{Y}^m, \mathcal{Q}^m)$, uniformly generate data points on the computational domain $\Omega_\text{I}$, resulting in a total of $20{,}000$ collocation points in $\mathcal{Y}^m$. For each $y_i^m$, we perform 500 Monte Carlo simulations to estimate $Q(y_i^m)$ using the stochastic representation \eqref{eq:sto_rep}, with parameters {$T_1 = 2$, $T_3 = 10$, and $\rho = 0.4$}. 
Then we sample $40{,}000$ points using Latin hypercube sampling (LHS) in the interior of $\Omega_\text{I}$ to minimize the residual of the Feynman--Kac equation. 
We use LHS to generate two additional training sets for the extension of $Q$ outside the domain $\Omega_\text{I}$. One contains $10{,}000$ external sample points in $\Omega_\text{O}^c$ and the corresponding constant extension $C=48$. The other is the intermediate set where we sample $10{,}000$ points in $\Omega_\text{O}\setminus\Omega_\text{I}$ to control $\mathcal{L}Q$.
We approximate the solution using a feedforward neural network with 6 hidden layers containing 16, 128, 128, 128, 32, and 8 neurons, respectively.
\begin{figure}[htp]
    \centering
    \includegraphics[width=0.27\linewidth]{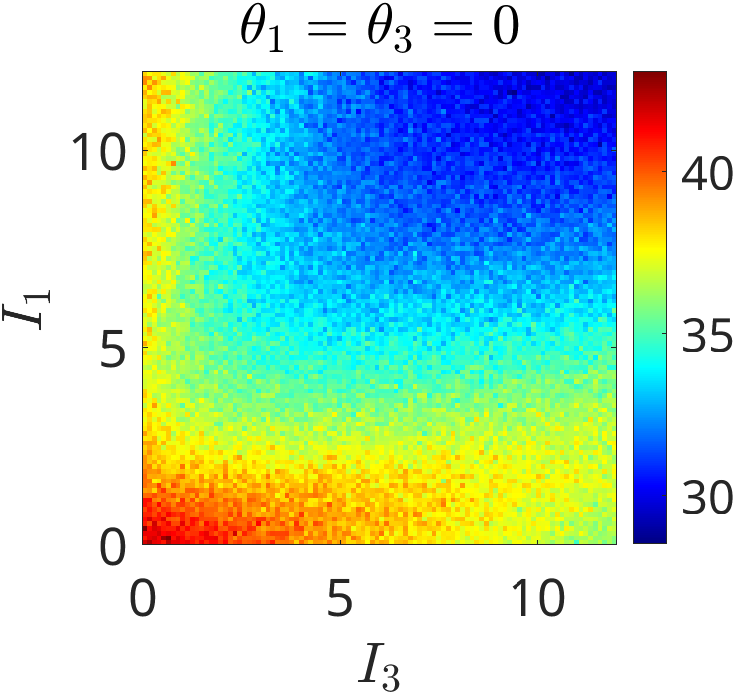}
    \includegraphics[width=0.27\linewidth]{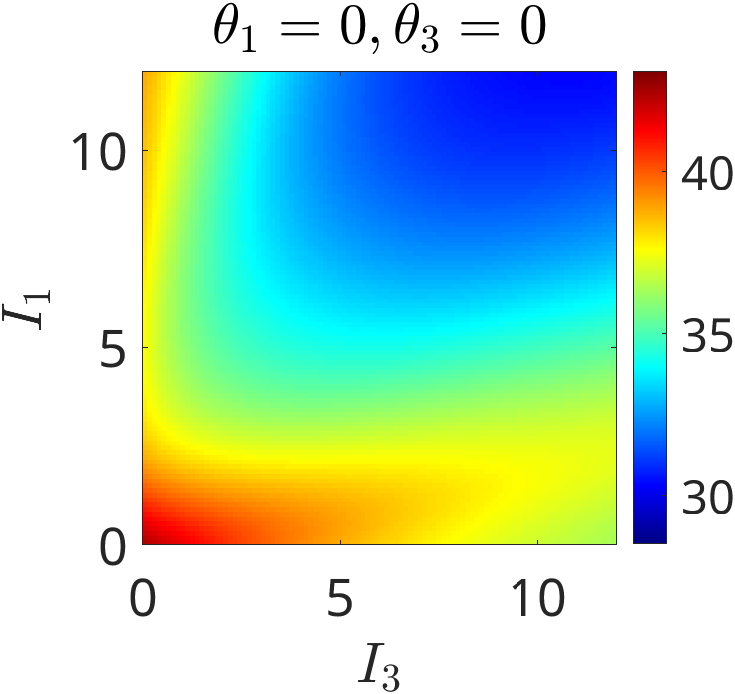}
    \caption{Computed $Q$ from Monte Carlo simulation and predicted solution from neural network, plotted on $(I_1, I_3)$-slice for $(\theta_1,\theta_3)=(0,0).$}
    \label{fig: I1I3slice}
\end{figure}
\begin{figure}[htp]
    \centering
    \includegraphics[width=0.24\linewidth]{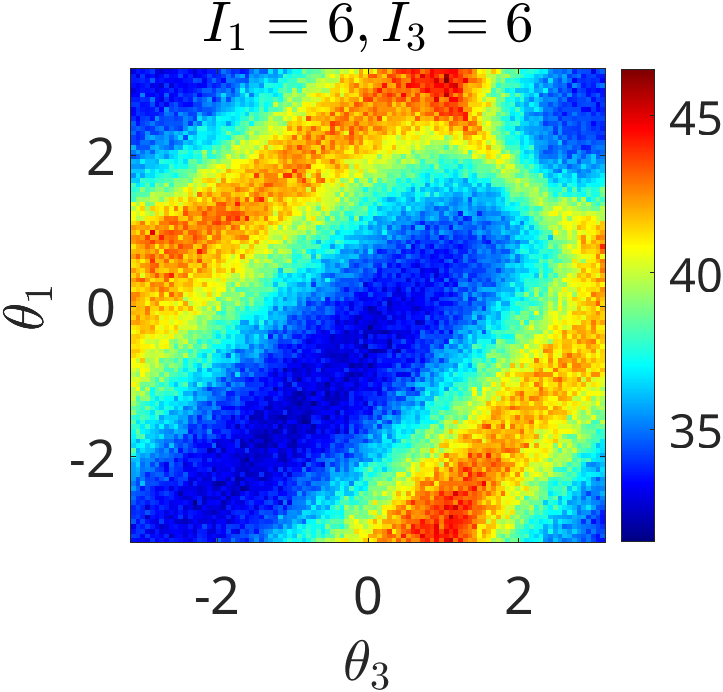}
    \includegraphics[width=0.24\linewidth]{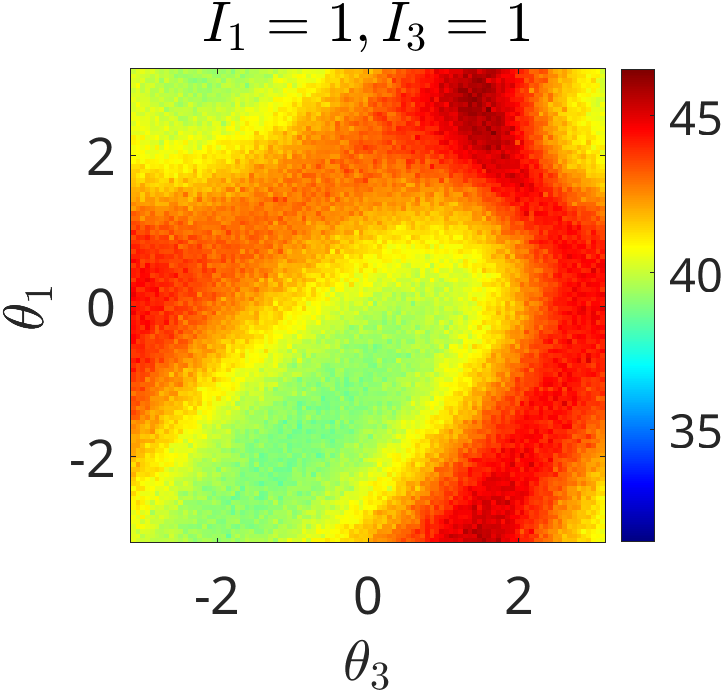} 
    \includegraphics[width=0.24\linewidth]{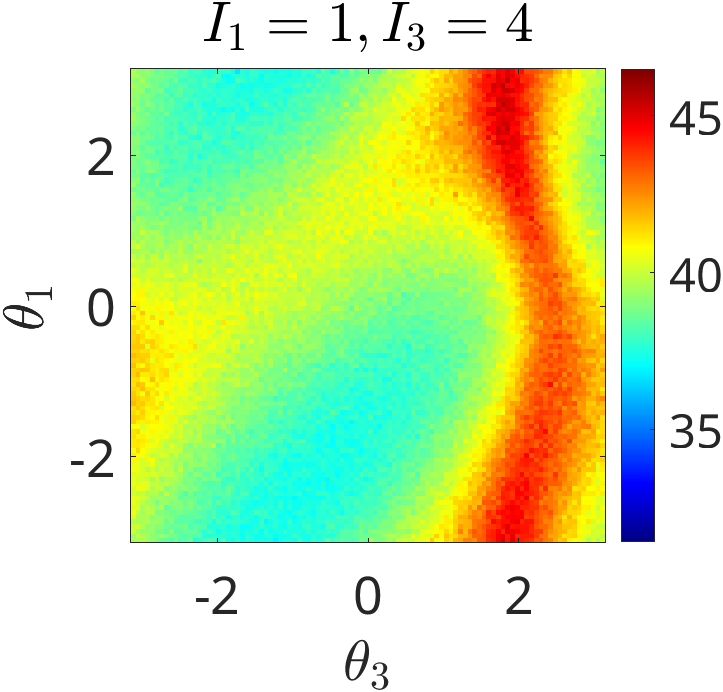}
    \includegraphics[width=0.24\linewidth]{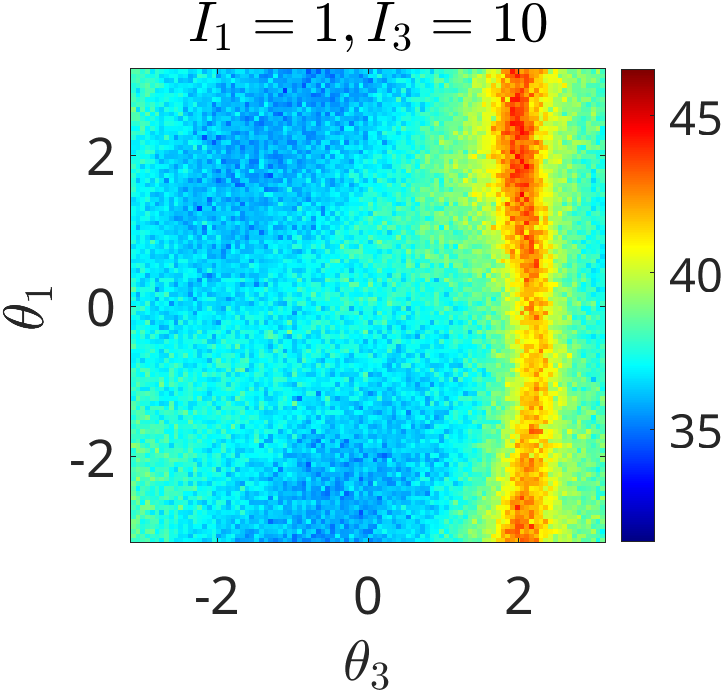}
    \caption{Computed $Q$ from Monte Carlo simulation, plotted on $(\theta_1, \theta_3)$-slices for different values of $I_1$ and $I_3$.}
    \label{fig: mc4d}
\end{figure}
\begin{figure}[htp]
    \centering
    \includegraphics[width=0.24\linewidth]{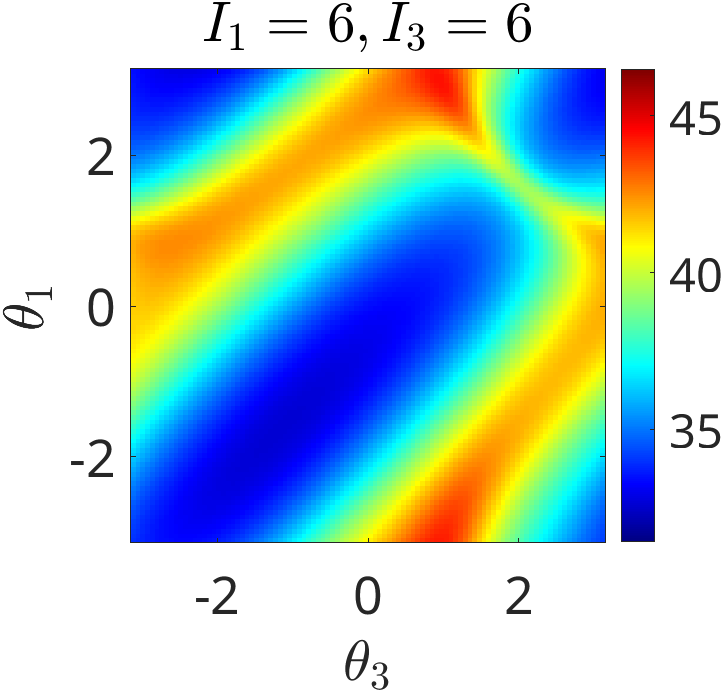}
    \includegraphics[width=0.24\linewidth]{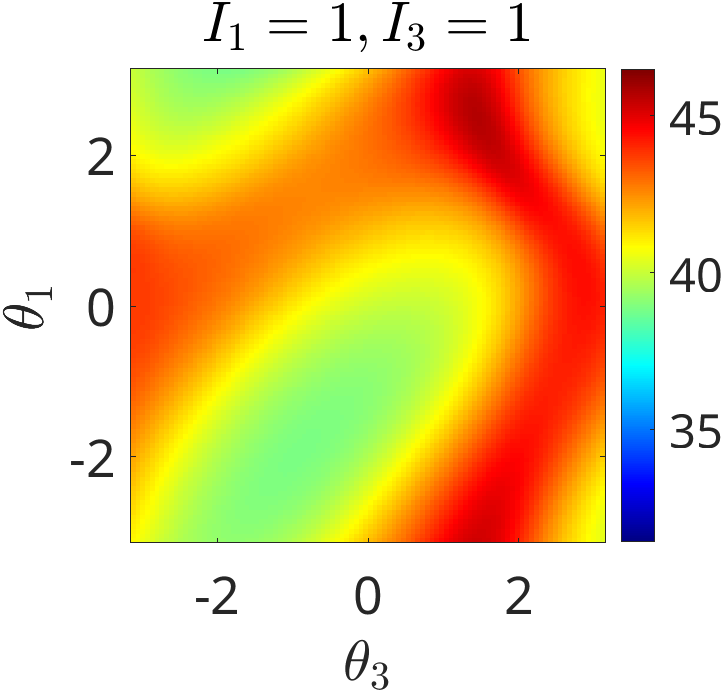}
    \includegraphics[width=0.24\linewidth]{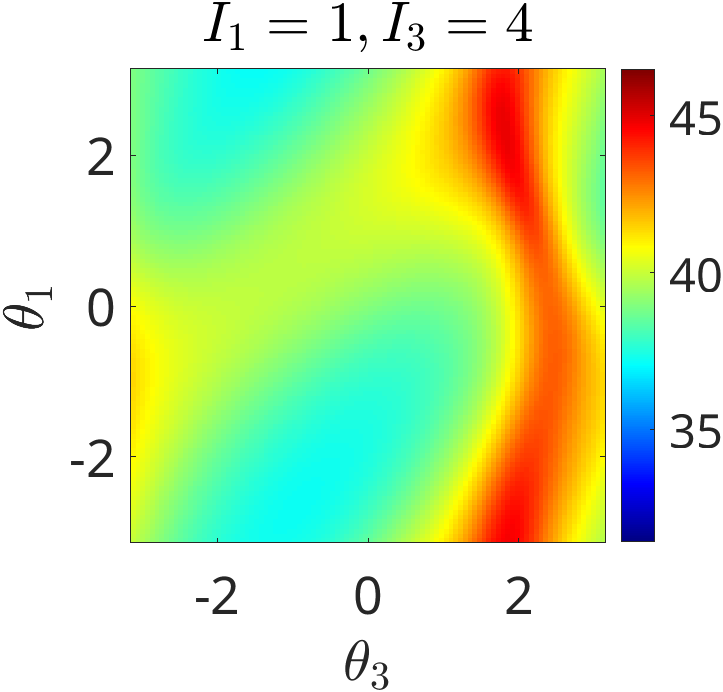}
    \includegraphics[width=0.24\linewidth]{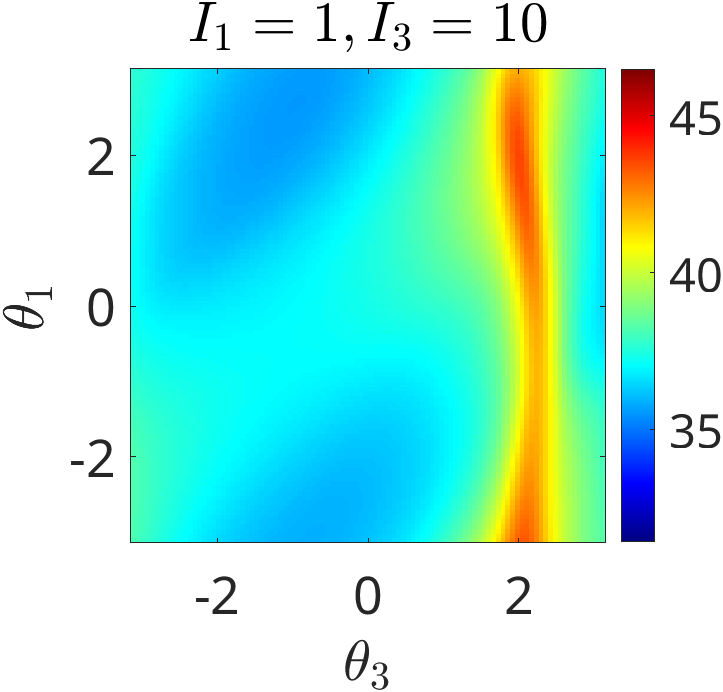}
    \caption{Numerical solution $Q$ learned from the neural network method, plotted on $(\theta_1, \theta_3)$-slices for different values of $I_1$ and $I_3$.}
    \label{fig: predicted}
\end{figure}

We compare the prefactor component $Q$ at a few two-dimensional slices of the domain at fixed $(\theta_1, \theta_3)$ in \Cref{fig: I1I3slice} and at fixed $(I_1, I_3)$ in \Cref{fig: mc4d} and \Cref{fig: predicted} with the Monte Carlo simulation result. \Cref{fig: mc4d} presents the results from high-resolution Monte Carlo simulations at $4$ selected $(\theta_1, \theta_3)$ slices and \Cref{fig: predicted} shows the corresponding learned function $Q$ from the neural network.
We observe that the neural network solution is consistent with the high-resolution Monte Carlo simulation, even though the Monte Carlo training data are relatively sparse. 
\begin{figure}[htp]
    \centering
    \includegraphics[width=0.26\linewidth]{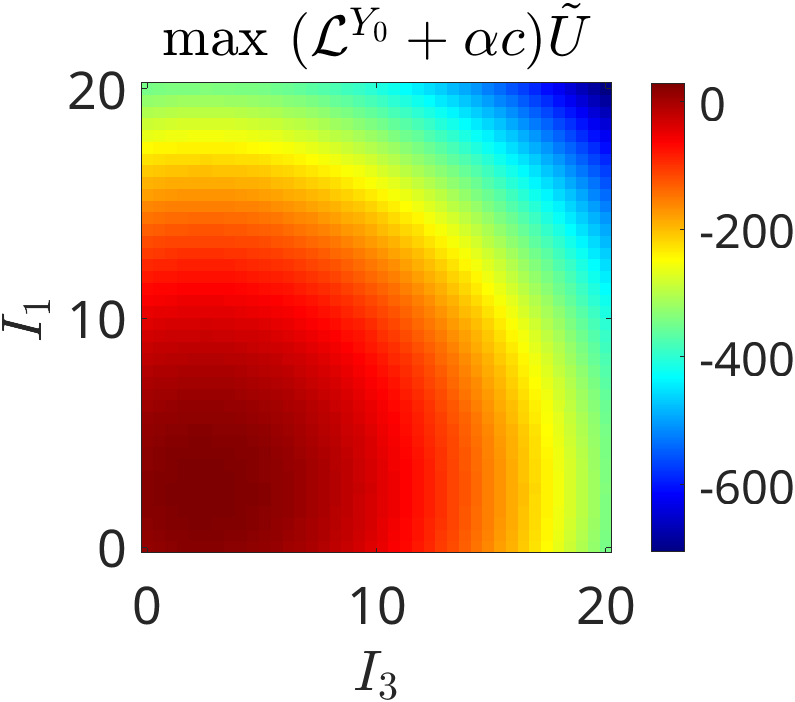}
    \includegraphics[width=0.26\linewidth]{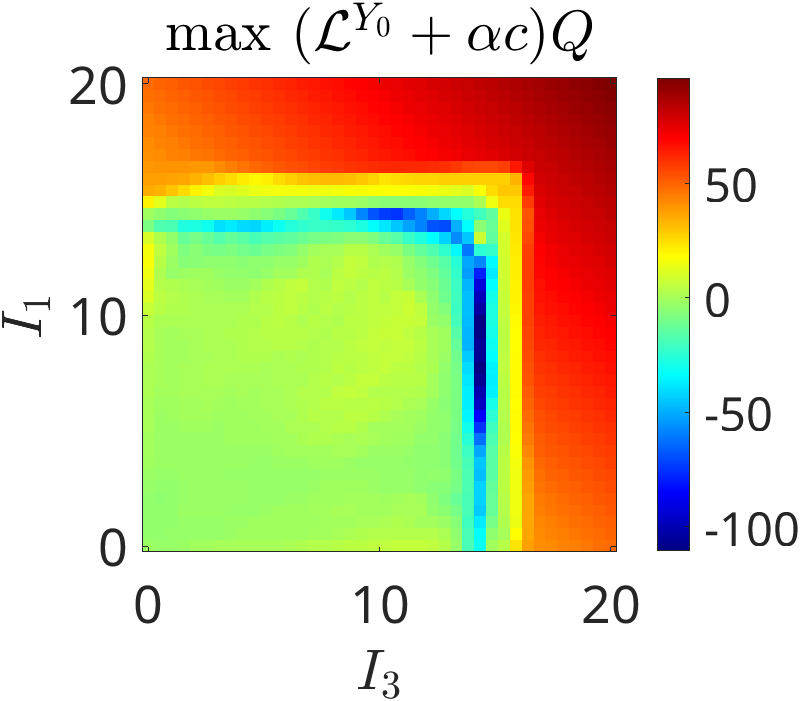}
    \includegraphics[width=0.26\linewidth]{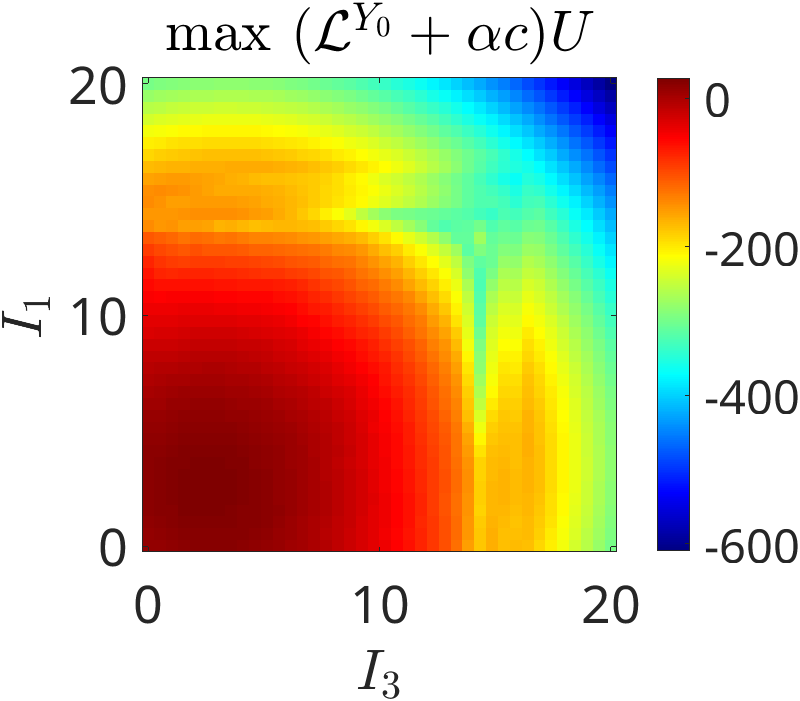}
    \caption{Panel Left, Middle, and Right: Maximum values of $(\mathcal{L}^{Y_0}+\alpha c)\tilde{U}$, $(\mathcal{L}^{Y_0}+\alpha c)Q$, and $(\mathcal{L}^{Y_0}+\alpha c){U}$ for $(I_1, I_3)$, with $I_1,I_3\in [0,20]$ at selected $(\theta_1, \theta_3)$ slices.}
    \label{fig: generator}
\end{figure}

Then we check \Cref{assumption1}(D), which generally can only be checked numerically. Since the extension of $Q$ is defined outside of $\Omega_\text{I}$, it is sufficient to check $\tilde{U}=3(I_1 + I_3)$ without the modification for small $I_1, I_3$. The left panel of \Cref{fig: generator} shows $\mathcal{L}^{Y_0} \tilde{U}$, which only depends on $(I_1, I_3)$. The middle panel of \Cref{fig: generator} shows the maximum of $\mathcal{L}^{Y_0}Q$ for the predicted $Q$ at each $(\theta_1, \theta_3)$ slice centered at $40 \times 40$ $(I_1, I_3)$ coordinates. We can see that $\mathcal{L}^{Y_0}Q$ is dominated by $\mathcal{L}^{Y_0} \tilde{U}$ both in $\Omega_\text{O}\setminus \Omega_\text{I}$ and in $\Omega_\text{O}^c$. As shown in the right panel of \Cref{fig: generator}, the value of $\mathcal{L}^{Y_0}U$, where $U=Q+\tilde{U}$, is uniformly negative in $\Omega_\text{I}^c$.
 
We can further relax the assumptions using the numerical solution of $Q$. We compute $U$ numerically in the bounded domain $\Omega_\text{I}$ and show that $\tilde{U}$ is a suitable exterior Lyapunov function. Recall that \Cref{assumption1}(E) can be relaxed to \Cref{assumption1}($\tilde{\text{E}}$) if $\|\nabla U\|_2$ is known.  For the three-mode energy cascade model, we have 
\begin{equation}
h= \| g - g_0\|_2 \leq \sqrt{288(1+\rho)^2 + 18}I_2\approx 47.7 I_2.
\end{equation}
Numerically, we compute $\|\nabla Q\|_{\infty} = 16.5035$, where the maximum is taken at 
\begin{equation}(I_1,I_3,\theta_1,\theta_3)=(9.2529,  11.4814,     2.09685,  1.8199).\end{equation}
 In addition, since $\tilde{U}$ is explicitly given, we easily have $\|\nabla \tilde{U}\|_2 = 13.9573$. Considering the worst case, we have an upper bound $\|\nabla U \|_2 = 30.4608$.
The strict bound in \Cref{assumption1}($\tilde{\text{E}}$) implies that $h(x,y)$ must be small enough such that $h \|\nabla U\|_2 \leq \gamma/2$. Therefore, $I_2^{-0.025}U$ is a Lyapunov function in the range $0 < I_2 < 4.52 \times 10^{-5}$.

\section{Conclusion}

The Feynman--Kac--Lyapunov method was first introduced in \cite{hani2025nonequilibrium} to address the stochastic stability problem of a specific energy cascade model. In this paper, we extend this method to a broader class of high-dimensional stochastic driving-driven systems, where the driven variables are approximately solely driven by the driving variables. The driving variables have a Lyapunov function that provides control over their stochastic stability, allowing us to construct a Lyapunov function for the full coupled system. 
The main idea is to construct a prefactor function that assigns larger values to the states of the driving dynamics that contribute to the instability of the driven dynamics. The product of this prefactor and the Lyapunov function of the driven variables then forms a Lyapunov function for the full system. This framework applies to a broad class of high-dimensional stochastic systems encountered in various applications. A Lyapunov function constructed from the infinitesimal generator offers advantages over the time-discrete drift condition obtained by averaging the driving dynamics. We provide general conditions under which such driving-driven systems admit a Lyapunov function.

If the system is one-way coupled, the prefactor is the eigenfunction of the Feynman--Kac operator. Otherwise, if the system is only an approximation of a one-way coupled system, one needs to solve the Feynman--Kac equation numerically and extend it into the full space. Therefore, we develop a numerical approach to solve the high-dimensional Feynman--Kac equation that determines the prefactor of the Lyapunov function. The method approximates the solution using a neural network, where the training process minimizes both the PDE residual and the discrepancy from the training data at preselected collocation points obtained from Monte Carlo simulations. The numerically computed prefactor function provides several advantages. It reveals which components of the driving dynamics contribute most to the instability of the driven system. Furthermore, once we solve the Feynman--Kac equation numerically and obtain the upper bound of the gradient norm, we can relax the theoretical conditions required for the proof. We demonstrate these properties through our numerical examples.

We can further extend the Feynman--Kac--Lyapunov method, which still has room for development in several directions.
The first possible extension concerns stabilization by noise. When the driven component $X_t$ possesses only a weak Lyapunov function $V$ satisfying $\nabla V \cdot \mathrm{d}X_t = 0$ on a submanifold $M$, our current result does not directly apply. However, we still expect $X_t$ to be stochastically stable, since the dynamics of $Y_t$ continuously drive $X_t$, preventing it from remaining on $M$ indefinitely. This setting corresponds to the phenomenon of stabilization by noise studied in \cite{bedrossian2024existence, bedrossian2024stationary, herzog2015noise1, herzog2015noise2}.
The second possible extension concerns the rigorous estimation of an upper bound for $|\nabla U|$ in the Feynman--Kac equation. Functional inequalities such as the Poincar\'e and Sobolev inequalities provide only lower bounds for $|\nabla U|$. A more effective approach is to use the stochastic representation of $\nabla U$ given by the Bismut--Elworthy--Li formula \cite{elworthy1994formulae}, which allows direct estimation of the gradient through probabilistic techniques. Our current approach needs to solve a Feynman--Kac equation on a bounded domain and then numerically calculate a $C^2$ extension to the whole domain because we do not have a good asymptotic control of the gradient of the eigenfunction of the Feynman--Kac operator as $|y| \rightarrow \infty$. The Bismut--Elworthy--Li formula can significantly help the gradient estimation, and potentially extend the applicability of a fully rigorous proof from a one-way coupled system like equation \eqref{SDE3} to some general cases like equation \eqref{SDE2}. We plan to address these extensions in our future work.

\section*{Acknowledgments}
W.C. was supported by NSF grants DMS-2350325 and DMS-2514053. Y.L. was partially supported by NSF grant DMS-2510209.

\bibliographystyle{plain}
\bibliography{reference} 

@article{nirenberg,
author = {Berestycki, H. and Nirenberg, L. and Varadhan, S. R. S.},
title = {The principal eigenvalue and maximum principle for second-order elliptic operators in general domains},
journal = {Communications on Pure and Applied Mathematics},
volume = {47},
number = {1},
pages = {47-92},
doi = {https://doi.org/10.1002/cpa.3160470105},
year = {1994}
}

@book{kato2013perturbation,
  title={Perturbation theory for linear operators},
  author = {Kato, Tosio},
  year={2013},
  publisher={Springer Berlin Heidelberg}
}

@book{Oksendal,
  title={Stochastic Differential Equations: An Introduction with Applications},
  author={Bernt Øksendal},
  series={Universitext},
  year={2003},
  publisher={Springer Berlin Heidelberg}
}

@book{meyn1993MC,
  title={{Markov} Chains and Stochastic Stability},
  author={Sean P. Meyn and Richard L. Tweedie},
  series={Communications and Control Engineering},
  year={1993},
  publisher={Springer London}
}

@article{hening2018coexistence,
  author = {Alexandru Hening and Dang H. Nguyen},
  title = {{Coexistence and extinction for stochastic Kolmogorov systems}},
  volume = {28},
  journal = {The Annals of Applied Probability},
  number = {3},
  publisher = {Institute of Mathematical Statistics},
  year = {2018},
  doi = {10.1214/17-AAP1347},
}

@article{liu2002representation,
    author = {Liu, Genqian},
    title = {Strongly Continuous Semigroups and Stochastic Representation},
    journal = {Journal of the London Mathematical Society},
    volume = {65},
    number = {3},
    pages = {639-660},
    year = {2002},
    month = {06},
    doi = {10.1112/S0024610701003003},
}

@book{evans2022partial,
  title={Partial Differential Equations},
  author={Evans, Lawrence C.},
  publisher={American Mathematical Society},
  year={2010},
  volume={19},
  series={Graduate Studies in Mathematics},
}

@article{meyn2003spectral,
  author = {I. Kontoyiannis and S. P. Meyn},
  title = {{Spectral theory and limit theorems for geometrically ergodic {Markov} processes}},
  volume = {13},
  journal = {The Annals of Applied Probability},
  number = {1},
  publisher = {Institute of Mathematical Statistics},
  year = {2003},
  doi = {10.1214/aoap/1042765670},
}

@misc{hani2025nonequilibrium,
      title={Non-equilibrium steady state for a three-mode energy cascade model}, 
      author={Zaher Hani and Yao Li and Andrea Nahmod and Gigliola Staffilani},
      year={2025},
      eprint={2505.16018},
      archivePrefix={arXiv},
      primaryClass={math.PR},
      url={https://arxiv.org/abs/2505.16018}, 
}

@InProceedings{zhai2022deep,
  title = {A deep learning method for solving {Fokker--Planck} equations},
  author =  {Zhai, Jiayu and Dobson, Matthew and Li, Yao},
  booktitle = {Proceedings of the 2nd Mathematical and Scientific Machine Learning Conference},
  pages = {568--597},
  year = 	{2022},
  volume = 	{145},
  series = 	{Proceedings of Machine Learning Research},
  publisher =  {PMLR},
}

@article{rosenthal1995minorization,
 author = {Jeffrey S. Rosenthal},
 journal = {Journal of the American Statistical Association},
 number = {430},
 pages = {558--566},
 publisher = {[American Statistical Association, Taylor & Francis, Ltd.]},
 title = {Minorization Conditions and Convergence Rates for {{Markov} Chain Monte Carlo}},
 volume = {90},
 year = {1995}
}

@article{meyn1994convergenceRates,
 author = {Sean P. Meyn and R. L. Tweedie},
 journal = {The Annals of Applied Probability},
 number = {4},
 pages = {981--1011},
 publisher = {Institute of Mathematical Statistics},
 title = {Computable Bounds for Geometric Convergence Rates of {{Markov}} Chains},
 volume = {4},
 year = {1994}
}

@book{khasminskii2012stochastic,
  title={Stochastic stability of differential equations},
  author={Khasminskii, Rafail},
  year={2012},
  publisher={Springer}
}

@article{hairer2006ergodicity,
  title={Ergodicity of the {2D Navier--Stokes} equations with degenerate stochastic forcing},
  author={Hairer, Martin and Mattingly, Jonathan C},
  journal={Annals of Mathematics},
  pages={993--1032},
  year={2006},
  publisher={JSTOR}
}

@article{herzog2015noise1,
  title={Noise-induced stabilization of planar flows {I}},
  author={Herzog, David P and Mattingly, Jonathan C},
  journal={Electron. J. Probab},
  volume={20},
  number={111},
  pages={1--43},
  year={2015}
}

@article{mattingly2002ergodicity,
  title={Ergodicity for {SDEs} and approximations: locally {Lipschitz} vector fields and degenerate noise},
  author={Mattingly, Jonathan C and Stuart, Andrew M and Higham, Desmond J},
  journal={Stochastic Processes and Their Applications},
  volume={101},
  number={2},
  pages={185--232},
  year={2002},
  publisher={Elsevier}
}

@book{stroock2008partial, 
  place={Cambridge}, 
  series={Cambridge Studies in Advanced Mathematics}, 
  title={Partial Differential Equations for Probabilists}, 
  publisher={Cambridge University Press}, 
  author={Stroock, Daniel W.}, 
  year={2008}, 
  collection={Cambridge Studies in Advanced Mathematics}
}

@article{herzog2015noise2,
  title={Noise-induced stabilization of planar flows {II}},
  author={Herzog, David P and Mattingly, Jonathan C},
  journal={Electron. J. Probab},
  volume={20},
  number={113},
  pages={1--37},
  year={2015}
}

@article{bedrossian2024stationary,
  title={Stationary measures for stochastic differential equations with degenerate damping},
  author={Bedrossian, Jacob and Liss, Kyle},
  journal={Probability Theory and Related Fields},
  volume={189},
  number={1},
  pages={101--178},
  year={2024},
  publisher={Springer}
}

@article{bedrossian2024existence,
  title={Existence of stationary measures for partially damped {SDEs} with generic, Euler-type nonlinearities},
  author={Bedrossian, Jacob and Blumenthal, Alex and Callis, Keagan and Liss, Kyle},
  journal={arXiv preprint arXiv:2407.16592},
  year={2024}
}

@article{baxendale1988large,
  title={Large deviations and stochastic flows of diffeomorphisms},
  author={Baxendale, Peter H and Stroock, Daniel W},
  journal={Probability Theory and Related Fields},
  volume={80},
  number={2},
  pages={169--215},
  year={1988},
  publisher={Springer}
}

@article{bedrossian2022almost,
  title={Almost-sure exponential mixing of passive scalars by the stochastic {Navier--Stokes} equations},
  author={Bedrossian, Jacob and Blumenthal, Alex and Punshon-Smith, Samuel},
  journal={The Annals of Probability},
  volume={50},
  number={1},
  pages={241--303},
  year={2022},
  publisher={Institute of Mathematical Statistics}
}

@article{bedrossian2025non,
  title={Non-uniqueness of stationary measures for stochastic systems with almost surely invariant manifolds},
  author={Bedrossian, Jacob and Blumenthal, Alex and Punshon-Smith, Sam},
  journal={arXiv preprint arXiv:2505.22903},
  year={2025}
}

@article{dolgopyat2004sample,
  title={Sample path properties of the stochastic flows},
  author={Dolgopyat, Dmitry and Kaloshin, Vadim and Koralov, Leonid},
  journal={The Annals of Probability},
  volume={32},
  number={1A},
  pages={1--27},
  year={2004},
  publisher={Institute of Mathematical Statistics}
}

@article{meng2020composite,
  title={A composite neural network that learns from multi-fidelity data: Application to function approximation and inverse PDE problems},
  author={Meng, Xuhui and Karniadakis, George Em},
  journal={Journal of Computational Physics},
  volume={401},
  pages={109020},
  year={2020},
  publisher={Elsevier}
}

@article{jagtap2020adaptive,
  title={Adaptive activation functions accelerate convergence in deep and physics-informed neural networks},
  author={Jagtap, Ameya D and Kawaguchi, Kenji and Karniadakis, George Em},
  journal={Journal of Computational Physics},
  volume={404},
  pages={109136},
  year={2020},
  publisher={Elsevier}
}

@article{wang2021understanding,
  title={Understanding and mitigating gradient flow pathologies in physics-informed neural networks},
  author={Wang, Sifan and Teng, Yujun and Perdikaris, Paris},
  journal={SIAM Journal on Scientific Computing},
  volume={43},
  number={5},
  pages={A3055--A3081},
  year={2021},
  publisher={SIAM}
}

@article{elworthy1994formulae,
  title={Formulae for the derivatives of heat semigroups},
  author={Elworthy, K David and Li, Xue-Mei},
  journal={Journal of Functional Analysis},
  volume={125},
  number={1},
  pages={252--286},
  year={1994},
  publisher={Elsevier}
}

@book{gilbarg1977elliptic,
  title={Elliptic partial differential equations of second order},
  author={Gilbarg, David and Trudinger, Neil S and Gilbarg, David and Trudinger, NS},
  volume={224},
  number={2},
  year={1977},
  publisher={Springer}
}

@article{lu2021physics,
  title={Physics-informed neural networks with hard constraints for inverse design},
  author={Lu, Lu and Pestourie, Raphael and Yao, Wenjie and Wang, Zhicheng and Verdugo, Francesc and Johnson, Steven G},
  journal={SIAM Journal on Scientific Computing},
  volume={43},
  number={6},
  pages={B1105--B1132},
  year={2021},
  publisher={SIAM}
}

@article{gao2021phygeonet,
  title={PhyGeoNet: Physics-informed geometry-adaptive convolutional neural networks for solving parameterized steady-state PDEs on irregular domain},
  author={Gao, Han and Sun, Luning and Wang, Jian-Xun},
  journal={Journal of Computational Physics},
  volume={428},
  pages={110079},
  year={2021},
  publisher={Elsevier}
}

@article{cai2021physics,
  title={Physics-informed neural networks (PINNs) for fluid mechanics: A review},
  author={Cai, Shengze and Mao, Zhiping and Wang, Zhicheng and Yin, Minglang and Karniadakis, George Em},
  journal={Acta Mechanica Sinica},
  volume={37},
  number={12},
  pages={1727--1738},
  year={2021},
  publisher={Springer}
}

@article{athreya2012propagating,
  title={Propagating {Lyapunov} functions to prove noise--induced stabilization},
  author={Athreya, Avanti and Kolba, Tiffany and Mattingly, Jonathan C},
  journal={Electron. J. Probab},
  volume={17},
  number={96},
  pages={1--38},
  year={2012}
}

@book{freidlin1985functional,
  title={Functional integration and partial differential equations},
  author={Freidlin, Mark Iosifovich},
  number={109},
  year={1985},
  publisher={Princeton University Press}
}

@article{wang2024deep,
  title={Deep adaptive sampling for surrogate modeling without labeled data},
  author={Wang, Xili and Tang, Kejun and Zhai, Jiayu and Wan, Xiaoliang and Yang, Chao},
  journal={Journal of Scientific Computing},
  volume={101},
  number={3},
  pages={77},
  year={2024},
  publisher={Springer}
}

@article{mSIR2024,
author = {Jiang, Ning and Chu, Weiqi and Li, Yao},
title = {Modeling, Inference, and Prediction in Mobility-Based Compartmental Models for Epidemiology},
journal = {SIAM Journal on Applied Mathematics},
volume = {85},
number = {5},
pages = {2355-2375},
year = {2025},
doi = {10.1137/24M1691557},
URL = { https://doi.org/10.1137/24M1691557},
eprint = { https://doi.org/10.1137/24M1691557},
}

@article{colliander2010transfer,
  title={Transfer of energy to high frequencies in the cubic defocusing nonlinear {S}chr{\"o}dinger equation},
  author={Colliander, James and Keel, Markus and Staffilani, Gigiola and Takaoka, Hideo and Tao, Terence},
  journal={Inventiones Mathematicae},
  volume={181},
  number={1},
  pages={39--113},
  year={2010},
  publisher={Springer}
}

\end{document}